\documentclass[a4paper,10pt]{article}
\usepackage{stmaryrd}
\usepackage{amsfonts}
\usepackage{bbm}
\usepackage{amscd}
\usepackage{mathrsfs}
\usepackage{latexsym,amssymb,amsmath,amscd,amscd,amsthm,amsxtra}
\usepackage[dvips]{graphicx}
\usepackage[utf8]{inputenc}
\usepackage[T1]{fontenc}
\usepackage{lmodern}
\usepackage{amssymb}
\usepackage[all]{xy}
\usepackage{nicefrac,mathtools,enumitem}
\usepackage{microtype}
\usepackage{xcolor}
\newtheorem{thm}{Theorem}[section]
\newtheorem{lem}[thm]{Lemma}
\newtheorem{cor}[thm]{Corollary}
\newtheorem{pro}[thm]{Proposition}
\newtheorem{ex}[thm]{Example}
\newtheorem{rmk}[thm]{Remark}
\newtheorem{defi}[thm]{Definition}

\newcommand {\emptycomment}[1]{} 

\newcommand{\be }{\begin{equation}}
\newcommand{\ee }{\end{equation}}

\newcommand{\pf}{\noindent{\bf Proof.}\ }

\newcommand{\Sym}{\mathsf{Sym}}

\newcommand{\Real}{\mathbb R}

\newcommand{\huaB}{\mathcal{B}}

\newcommand{\huaA}{\mathcal{A}}
\newcommand{\huaL}{\mathcal{L}}

\newcommand{\huaF}{\mathcal{F}}

\newcommand{\huaH}{\mathcal{H}}

\newcommand{\CWM}{C^{\infty}(M)}

\newcommand{\frkg}{\mathfrak g}

\newcommand{\frkB}{\mathfrak B}

\newcommand{\frkD}{\mathfrak D}

\newcommand{\frkX}{\mathfrak X}

\def\qed{\hfill ~\vrule height6pt width6pt depth0pt}

\newcommand{\half}{\frac{1}{2}}

\newcommand{\br}[1]{   [ \cdot,    \cdot  ]   }

\newcommand{\dev}{\mathfrak{D}}

\newcommand{\id}{\rm{id}}

\newcommand{\dM}{\mathrm{d}}

\newcommand{\Hom}{\mathrm{Hom}}

\newcommand{\End}{\mathrm{End}}

\newcommand{\Nat}{\mathbb N}

\newcommand{\K}{\mathbb{K}}

\newcommand{\KVN}{\mathrm{KVN}}

\newcommand{\KVB}{\mathrm{KV\Omega}}

\newcommand{\HN}{\mathrm{H N}}

\newcommand{\POm}{\mathrm{P\Omega}}

\newcommand {\yh}[1]{{\marginpar{*}\scriptsize\textcolor{red}{  #1}}}

\begin{document}
\title{Koszul-Vinberg structures and compatible structures on left-symmetric algebroids
\thanks
 {
This research  is supported by
NSF of China (11471139) and NSF of Jilin Province (20170101050JC).
 }}
\author{Qi Wang$^{1}$, Jiefeng Liu$^2$ and Yunhe Sheng$^{1}$\\
 $^1$Department of Mathematics, Jilin University,\\
 \vspace{2mm}Changchun 130012, Jilin, China\\
$^2$Department of Mathematics, Xinyang Normal University,\\
\vspace{2mm} Xinyang 464000, Henan, China\\
Email:wangqi17@mails.jlu.edu.cn; liujf12@126.com; shengyh@jlu.edu.cn\\
}

\date{}
\footnotetext{{\it{Keywords}: left-symmetric algebroid, Koszul-Vinberg structure, Koszul-Vinberg-Nijenhuis structure,  pseudo-Hessian-Nijenhuis structure }}
\footnotetext{{\it{MSC}}: 17B62, 53D17}

\maketitle
\begin{abstract}
 In this paper, we introduce the notion of  Koszul-Vinberg-Nijenhuis structures on a left-symmetric algebroid as analogues of   Poisson-Nijenhuis structures on a Lie algebroid, and show that a Koszul-Vinberg-Nijenhuis structure  gives rise to a hierarchy of Koszul-Vinberg structures. We introduce the notions of $\KVB$-structures, pseudo-Hessian-Nijenhuis structures and complementary symmetric $2$-tensors for Koszul-Vinberg structures on left-symmetric algebroids, which are  analogues of $\POm$-structures, symplectic-Nijenhuis structures and complementary 2-forms for Poisson structures. We also study the relationships between these various structures.
\end{abstract}

\tableofcontents
\section{Introduction}
Left-symmetric algebras (or pre-Lie algebras) arose from the
study of convex homogeneous cones (\cite{Vinb}), affine manifolds
and affine structures on Lie groups (\cite{Koszul1}), deformation
and cohomology theory of associative algebras (\cite{G}) and then
appear in many fields in mathematics and mathematical physics. See
the survey article \cite{Pre-lie algebra in geometry} and the
references therein. In particular, there are
 close relations between left-symmetric algebras and
certain important left-invariant structures on Lie groups
(\cite{Kim,Lichnerowicz,MT,Med,Milnor}).

A left-symmetric algebroid is a geometric generalization of a
left-symmetric algebra. See \cite{LiuShengBaiChen,Boyom1,Boyom2} for more details and applications. The notion of a Nijenhuis operator on a left-symmetric algebroid was introduced in \cite{lsb}, which could generate a trivial deformation. See \cite{WBLS}  for more details of deformations of left-symmetric algebras.
 In \cite{lsb2}, Motivated by the theory of Lie biagebroids (\cite{Lie bialgebroid}),   the notion of a left-symmetric bialgebroid was introduced as a geometric generalization of a left-symmetric bialgebra (\cite{Left-symmetric bialgebras}). The double of a left-symmetric bialgebroid is not a left-symmetric algebroid anymore, but a pre-symplectic algebroid (\cite{lsb}). This result is  parallel to the fact that
 the double of a Lie bialgebroid is a Courant algebroid (\cite{lwx}). As a Poisson structure $\pi$ on a manifold gives rise to a Lie bialgebroid,  a symmetric $2$-tensor $H$ on a flat manifold satisfying a certain condition  gives rise to a left-symmetric bialgebroid. In this paper, we call this symmetric $2$-tensor $H$ a Koszul-Vinberg structure. In particular,  if the Koszul-Vinberg structure $H$ is nondegenerate, the inverse of $H$ is  a pseudo-Hessian structure on a flat manifold. See \cite{NiBai,Shima,Geometry of Hessian structures} for more information about pseudo-Hessian Lie algebras and Hessian geometry. Therefore,   Koszul-Vinberg structures and pseudo-Hessian structures are symmetric analogues of  Poisson structures and symplectic structures, respectively. In \cite{BBo}, the authors recovered the Koszul-Vinberg structure (they called it a contravariant pseudo-Hessian structure) from a different point of view and built an analogue of Darboux-Weinstein theorem near a regular point for this structure.


In this paper, we add some compatibility conditions between a Koszul-Vinberg structure and a Nijenhuis operator on a left-symmetric algebroid to introduce the notion of a Koszul-Vinberg-Nijenhuis structure. Koszul-Vinberg-Nijenhuis structures on left-symmetric algebroids enjoy many properties that parallel to  Poisson-Nijenhuis structures on Lie algebroids (\cite{Kosmann1,Kosmann2,MM}). Furthermore, we introduce the notions of $\KVB$-structures, pseudo-Hessian-Nijenhuis structures and complementary symmetric $2$-tensors for Koszul-Vinberg structures on left-symmetric algebroids, which are analogues of $\POm$-structures (\cite{Kosmann2,MM}), symplectic-Nijenhuis structures
(\cite{Dorfman1993,Kosmann2,MM}) and complementary 2-forms for Poisson structures (\cite{Kosmann2,Vaisman}).

The paper is organized as follows. In Section \ref{sec:Prelimainaries}, we recall the notions and properties of Lie algebroids and left-symmetric algebroids. 
In Section \ref{sec:KVN}, we add some compatibility conditions between a Koszul-Vinberg structure $H$ and a Nijenhuis operator $N$ to introduce the notion of a Koszul-Vinberg-Nijenhuis structure ($\KVN$-structure) on a left-symmetric algebroid.
In Section \ref{sec:comgRB}, we study compatible Koszul-Vinberg structures on left-symmetric algebroids. We show that a $\KVN$-structure $(H,N)$ gives rise to a hierarchy of Koszul-Vinberg structures. 
In Section \ref{sec:MC},  we introduce the notion of a $\KVB$-structure on a left-symmetric algebroid, which consists of a Koszul-Vinberg structure $H$ and a    $2$-cocycle $\frkB\in\Sym^2(A^*)$ satisfying some compatibility conditions. 
 The relations between $\KVB$-structures and $\KVN$-structures are studied.
In Section \ref{sec:HN}, we introduce the notion of a pseudo-Hessian-Nijenhuis structure ($\HN$-structure) on a left-symmetric algebroid $A$, which consists of a pseudo-Hessian structure $\frkB\in\Sym^2(A^*)$ and a Nijenhuis operator $N$ on $A$ satisfying some compatibility conditions. The relations among $\KVN$-structures, $\KVB$-structures and $\HN$-structures are discussed. Moreover, various examples are given.

\section{Preliminaries}\label{sec:Prelimainaries}


\subsection{Lie algebroids}

The notion of a Lie algebroid was introduced
by Pradines in 1967, as a generalization of Lie algebras and
tangent bundles. See \cite{General theory of Lie groupoid and Lie
algebroid} for the general theory about Lie algebroids.
\begin{defi}
A {\bf Lie algebroid} structure on a vector bundle $\huaA\longrightarrow M$ is
a pair that consists of a Lie algebra structure $[-,-]_\huaA$ on
the section space $\Gamma(\huaA)$ and a  bundle map
$a_\huaA:\huaA\longrightarrow TM$, called the anchor, such that the
following relation is satisfied:
$$~[x,fy]_\huaA=f[x,y]_\huaA+a_\huaA(x)(f)y,\quad \forall ~x,y\in\Gamma(\huaA),~f\in
\CWM.$$
\end{defi}

For a vector bundle $E\longrightarrow M$, we denote by
$\dev(E)$   the gauge Lie algebroid of the
 frame bundle
 $\huaF(E)$, which is also called the covariant differential operator bundle of $E$.

Let $(\huaA,[-,-]_\huaA,a_\huaA)$ and $(\huaB,[-,-]_\huaB,a_\huaB)$ be two Lie
algebroids (with the same base), a {\bf base-preserving morphism}
from $\huaA$ to $\huaB$ is a bundle map $\sigma:\huaA\longrightarrow \huaB$ over the identity such
that
\begin{eqnarray*}
  a_\huaB\circ\sigma=a_\huaA,\quad
  \sigma[x,y]_\huaA=[\sigma(x),\sigma(y)]_\huaB.
\end{eqnarray*}

A {\bf representation} of a Lie algebroid $\huaA$
 on a vector bundle $E$ is a base-preserving morphism $\rho$ from $\huaA$ to the Lie algebroid $\dev(E)$.
Let $(E;\rho)$ be a representation.
The {\bf dual representation} of a Lie algebroid $\huaA$ on $E^*$ is the bundle map $\rho^*:\huaA\longrightarrow \dev(E^*)$ given by
$$
\langle \rho^*(x)(\alpha),u\rangle=a_\huaA(x)\langle \alpha,u\rangle-\langle \alpha,\rho(x)(u)\rangle,\quad \forall x\in \Gamma(\huaA),~\alpha\in\Gamma(E^*),u\in\Gamma(E).
$$

For all $x\in \Gamma(\huaA)$, the {\bf Lie derivation} $\huaL_x:\Gamma(\huaA^*)\longrightarrow\Gamma(\huaA^*)$ of the Lie algebroid $\huaA$  is given by
\begin{eqnarray}\label{Lie der1}
\langle\huaL_x\alpha,y\rangle=a_\huaA(x)\langle \alpha,y\rangle-\langle \alpha,[x,y]_\huaA\rangle,\quad \forall y\in\Gamma(\huaA),\alpha\in\Gamma(\huaA^*).
\end{eqnarray}
A Lie algebroid $(\huaA,[-,-]_\huaA,a_\huaA)$   naturally
represents on the trivial line bundle $E=M\times \mathbb R$ via
the anchor map $a_\huaA:\huaA\longrightarrow TM$. The
 coboundary operator
${\dM_\huaA}:\Gamma(\wedge^k\huaA^*)\longrightarrow
\Gamma(\wedge^{k+1}\huaA^*)$ is given by
\begin{eqnarray*}
  {\dM_\huaA}\varpi(x_1,\cdots,x_{k+1})&=&\sum_{i=1}^{k+1}(-1)^{i+1}a_\huaA(x_i)\varpi(x_1\cdots,\widehat{x_i},\cdots,x_{k+1})\\
  &&+\sum_{i<j}(-1)^{i+j}\varpi([x_i,x_j]_\huaA,x_1\cdots,\widehat{x_i},\cdots,\widehat{x_j},\cdots,x_{k+1}).
\end{eqnarray*}
\subsection{Left-symmetric algebroids}

\begin{defi}
A {\bf left-symmetric algebra} is a pair $(\frkg,\cdot_\frkg)$, where ${\frkg}$ is a vector space, and  $\cdot_\frkg:{\frkg}\otimes {\frkg}\longrightarrow{\frkg}$ is a bilinear operation
such that for all $x,y,z\in {\frkg}$, the associator
\begin{equation}\label{eq:associator}
(x,y,z)\triangleq x\cdot_\frkg(y\cdot_\frkg z)-(x\cdot_\frkg y)\cdot_\frkg z
\end{equation} is symmetric in $x,y$,
i.e.
$$(x,y,z)=(y,x,z),\;\;{\rm or}\;\;{\rm
equivalently,}\;\;x\cdot_\frkg(y\cdot_\frkg z)-(x\cdot_\frkg y)\cdot_\frkg z=y\cdot_\frkg(x\cdot_\frkg z)-(y\cdot_\frkg x)\cdot_\frkg
z.$$
\end{defi}

\begin{defi}{\rm(\cite{LiuShengBaiChen,Boyom1})}\label{defi:left-symmetric algebroid}
A {\bf left-symmetric algebroid} structure on a vector bundle
$A\longrightarrow M$ is a pair that consists of a left-symmetric
algebra structure $\cdot_A$ on the section space $\Gamma(A)$ and a
vector bundle morphism $a_A:A\longrightarrow TM$, called the anchor,
such that for all $f\in\CWM$ and $x,y\in\Gamma(A)$, the following
conditions are satisfied:
\begin{itemize}
\item[\rm(i)]$~x\cdot_A(fy)=f(x\cdot_A y)+a_A(x)(f)y,$
\item[\rm(ii)] $(fx)\cdot_A y=f(x\cdot_A y).$
\end{itemize}
\end{defi}

We usually denote a left-symmetric algebroid by $(A,\cdot_A, a_A)$.
Any left-symmetric  algebra is a left-symmetric algebroid over a point.

For any $x\in\Gamma(A)$, we define
$L_x:\Gamma(A)\longrightarrow\Gamma(A)$  and $R_x:\Gamma(A)\longrightarrow\Gamma(A)$ by
\begin{equation}\label{eq:leftmul}L_xy=x\cdot_A y,\quad R_xy=y\cdot_A x.\end{equation}
Condition (i) in the above definition means that $L_x\in \frkD(A)$.
Condition (ii) means that the map $x\longmapsto L_x$ is
$C^\infty(M)$-linear. Thus, $L:A\longrightarrow \frkD(A)$ is a bundle
map. With the same notations, there is two maps $L_x,R_x:\Gamma(A^*)\longrightarrow\Gamma(A^*)$ given by
\begin{equation}\label{eq:dualLR}
\langle L_x \xi,y\rangle=a_A(x)\langle\xi,y\rangle-\langle \xi,L_x y\rangle,\quad \langle R_x \xi,y\rangle=-\langle \xi,R_x y\rangle,\quad \forall x,y\in\Gamma(A), \xi\in\Gamma(A^*).
\end{equation}
\begin{pro}{\rm (\cite{LiuShengBaiChen})}\label{thm:sub-adjacent}
  Let $(A,\cdot_A, a_A)$ be a left-symmetric algebroid. Define  a skew-symmetric bilinear bracket operation $[-,-]_A$ on $\Gamma(A)$ by
  $$
  [x,y]_A=x\cdot_A y-y\cdot_A x,\quad \forall ~x,y\in\Gamma(A).
  $$
Then $(A,[-,-]_A,a_A)$ is a Lie algebroid, and denoted by
$A^c$, called the {\bf sub-adjacent Lie algebroid} of
 $(A,\cdot_A,a_A)$. Furthermore, $L:A\longrightarrow \frkD(A)$  gives a
  representation of the Lie algebroid  $A^c$.
\end{pro}
Recall that a connection $\nabla$ on $M$ is said to be {\bf flat} if the torsion tensor and curvature tensor of $\nabla$ vanish identically. A manifold $M$ endowed with a flat connection $\nabla$ is called a {\bf flat manifold}.

The following results for flat manifolds are well known. See \cite{Shima} for more details.
\begin{pro}
  Let $(M,\nabla)$ be a flat manifold. Then there exists local coordinate system $\{x^1,\cdots,x^n\}$ on $M$ such that
  $\nabla_{\frac{\partial}{\partial x^i}}{\frac{\partial}{\partial
x^j}}=0,~ i,j=1,\ldots,n.$
We call this local coordinate system an {\bf affine coordinate system} with respect to $\nabla$.
\end{pro}

\begin{ex}
 Let $\Real^n$ be the affine space with the affine coordinate system $\{x^1,x^2,\ldots,x^n\}$. Then there is a natural flat connection $\nabla$ on $\Real^n$ defined by
 $\nabla_{\frac{\partial}{\partial x^i}}{\frac{\partial}{\partial
x^j}}=0,~i,j=1,\ldots,n.$
 We call this connection the {\bf standard flat connection} on $\Real^n$.
 \end{ex}

\begin{ex}\label{ex:main}{\rm Let $(M,\nabla)$ be a flat manifold. Then $(TM,\nabla,\id)$ is a left-symmetric algebroid whose sub-adjacent Lie algebroid is
exactly the tangent Lie algebroid. We denote this left-symmetric algebroid by $T_\nabla M$.
}
\end{ex}

In the following, we recall the definition of a Nijenhuis operator on a left-symmetric algebroid, which could generate a trivial deformation of the left-symmetric algebroid.

\begin{defi}{\rm(\cite{LiuShengBaiChen})}
A bundle map
$N:A\longrightarrow A$ is called a {\bf Nijenhuis operator} on a left-symmetric algebroid $(A,\cdot_A,a_A)$  if the Nijenhuis condition holds:
\begin{eqnarray}\label{integral condition of Nij}
N(x)\cdot_A N(y)-N(x\cdot_A N(y)+N(x)\cdot_A y+N(x\cdot_A
y))=0, \quad \forall x,y\in \Gamma(A).
\end{eqnarray}

\end{defi}
Obviously, any Nijenhuis operator on a left-symmetric algebroid is also a Nijenhuis operator on the corresponding  sub-adjacent Lie algebroid.

Let $N$ be a Nijenhuis operator on a left-symmetric algebroid $({A},\cdot_{A},a_A)$.  We denote by $\cdot_N:{A}\otimes {A}\longrightarrow{A}$ the deformed multiplication. More precisely,
\begin{equation}\label{eq:BNS3}
  x\cdot_N y=N(x)\cdot_{A} y+x\cdot_{A} N(y)-N(x\cdot_{A} y).
\end{equation}
Then $({A},\cdot_N,a_N=a_A\circ N)$ is a left-symmetric algebroid, and $N$ is a left-symmetric algebroid morphism between $({A},\cdot_N,a_N=a_A\circ N)$ and $(A,\cdot_A,a_A)$. We denote by $[-,-]_N$ the commutator bracket of $\cdot_N$. More precisely,
\begin{equation}\label{eq:deform-bracket}
  [x,y]_N=x\cdot_N y-y\cdot_N x,\quad\forall~x,y\in{\Gamma(A)}.
\end{equation}
Then $(A,[-,-]_N,a_N)$ is a Lie algebroid and $N$ is a Lie algebroid morphism between $({A},[-,-]_N,a_N=a_A\circ N)$ and $(A,[-,-]_A,a_A)$.

By direct calculations, we have
\begin{lem}\label{lem:Niejproperty}
  Let $({A},\cdot_{A},a_A)$ be a left-symmetric algebroid and $N$ a Nijenhuis operator.
  \begin{itemize}
\item[$\rm(i)$] For all $k\in\Nat$, $(A,\cdot_{N^k},a_{N^k}=a_A\circ N^k)$ is a left-symmetric algebroid;
\item[$\rm(ii)$]For all $l\in\Nat$, $N^l$ is  a Nijenhuis operator on the left-symmetric algebroid $(A,\cdot_{N^k},a_{N^k})$;
\item[$\rm(iii)$]The left-symmetric algebroids $(A,(\cdot_{N^k})_{N^l},a_{N^{k+l}})$ and $(A,\cdot_{N^{k+l}},a_{N^{k+l}})$ are the same;
\item[$\rm(iv)$]$N^l$ is a left-symmetric algebroid homomorphism from $(A,\cdot_{N^{k+l}},a_{N^{k+l}})$ to $(A,\cdot_{N^k},a_{N^{k}})$.
  \end{itemize}
\end{lem}

Let $(A,\cdot_A,a_A)$ be a left-symmetric algebroid and $E$  a vector
bundle. A {\bf representation} of $A$ on $E$ consists of a pair
$(\rho,\mu)$, where $\rho:A\longrightarrow \frkD(E)$ is a representation
of $A^c$ on $E $ and $\mu:A\longrightarrow \End(E)$ is a bundle
map, such that for all $x,y\in \Gamma(A),\ e\in\Gamma(E)$, we have
\begin{eqnarray}\label{representation condition 2}
 \rho(x)\mu(y)e-\mu(y)\rho(x)e=\mu(x\cdot_A y)e-\mu(y)\mu(x)e.
\end{eqnarray}
Denote a representation by $(E;\rho,\mu)$.

\emptycomment{let us recall the cohomology complex for a left-symmetric algebroid $(A,\cdot_A,a_A)$ with a representation $(E;\rho,\mu)$ briefly. Denote the set of $(n+1)$-cochains by
$$C^{n+1}(A,E)=\Gamma(\Hom(\wedge^{n}A\otimes A,E)),\
n{A}eq 0.$$  For all $\varphi\in C^{n}(A,E)$, the coboundary operator $\delta:C^{n}(A,E)\longrightarrow C^{n+1}(A,E)$ is given by
 \begin{eqnarray*}
\delta\varphi(x_1,\cdots,x_{n+1})&=&\sum_{i=1}^{n}(-1)^{i+1}\rho(x_i)\omega(x_1,\cdots,\hat{x_i},\cdots,x_{n+1})\\
 &&+\sum_{i=1}^{n}(-1)^{i+1}\mu(x_{n+1})\omega(x_1,\cdots,\hat{x_i},\cdots,x_n,x_i)\\
 &&-\sum_{i=1}^{n}(-1)^{i+1}\varphi(x_1,\cdots,\hat{x_i},\cdots,x_n,x_i\cdot_A x_{n+1})\\
 &&+\sum_{1\leq i<j\leq n}(-1)^{i+j}\varphi([x_i,x_j]_A,x_1,\cdots,\hat{x_i},\cdots,\hat{x_j},\cdots,x_{n+1}),
\end{eqnarray*}
for all $x_i\in \Gamma(A),i=1,\cdots,n+1$.
}

Let us recall the cochain complex with the coefficient in the trivial representation, i.e. $\rho=a_A$ and $\mu=0$. See \cite{cohomology of pre-Lie,LiuShengBaiChen} for general theory of cohomologies of right-symmetric algebras and left-symmetric algebroids respectively. The set of $(n+1)$-cochains  is given by
$$C^{n+1}(A)=\Gamma(\wedge^{n}A^*\otimes A^*),\
n\geq0.$$  For all $\varphi\in C^{n}(A)$ and $x_i\in
\Gamma(A),~i=1,\cdots,n+1$, the  coboundary operator
$\delta_A$ is given by
 \begin{eqnarray}\label{LSCA cohomology}
\nonumber\delta_A\varphi(x_1, \cdots,x_{n+1})
 &=&\sum_{i=1}^{n}(-1)^{i+1}a_A(x_i)\varphi(x_1, \cdots,\hat{x_i},\cdots,x_{n+1})\nonumber\\
 &&-\sum_{i=1}^{n}(-1)^{i+1}\varphi(x_1, \cdots,\hat{x_i},\cdots,x_n,x_i\cdot_A x_{n+1})\nonumber\\
 &&+\sum_{1\leq i<j\leq n}(-1)^{i+j}\varphi([x_i,x_j]_A,x_1,\cdots,\hat{x_i},\cdots,\hat{x_j},\cdots,x_{n+1}).
\end{eqnarray}

Let $(A,\cdot_A,a_A)$ be a left-symmetric algebroid. Define
$$\Sym^2(A)=\{H\in A\otimes A| H(\alpha,\beta)=H(\beta,\alpha),\quad \forall\alpha,\beta\in \Gamma(A^*)\}.$$
For any $H\in \Sym^2(A)$, the bundle map  $H^\sharp:A^*\longrightarrow A$  is given by $H^\sharp(\alpha)(\beta)=H(\alpha,\beta)$.
 We introduce   $\llbracket H,H\rrbracket_A\in\Gamma(\wedge^2A \otimes A) $ as follows:
\begin{eqnarray}\label{brac2}
&&\llbracket H,H\rrbracket_A(\alpha_1,\alpha_2,\alpha_3)=a_A(H^\sharp(\alpha_1))\langle H^\sharp(\alpha_2),\alpha_3\rangle-a_A(H^\sharp(\alpha_2))\langle H^\sharp(\alpha_1),\alpha_3\rangle\nonumber\\
&&\quad\quad\quad+\langle \alpha_1,H^\sharp(\alpha_2)\cdot_A H^\sharp(\alpha_3)\rangle-\langle\alpha_2,H^\sharp(\alpha_1)\cdot_A H^\sharp(\alpha_3)\rangle-\langle \alpha_3,[H^\sharp(\alpha_1),
H^\sharp(\alpha_2)]_A\rangle.
\end{eqnarray}
 Suppose that $H^\sharp:A^*\longrightarrow A$ is   nondegenerate. Then   $(H^\sharp)^{-1}:A\longrightarrow A^*$ is also a symmetric bundle map, which gives rise to an element, denoted by $H^{-1}$,  in $\Sym^2(A^*)$.
\begin{pro}{\rm(\cite{lsb2})}\label{pro:equivelent}
Let $(A,\cdot_A,a_A)$ be a left-symmetric algebroid and $H\in \Sym^2(A)$. If $H$ is nondegenerate, then
 $\llbracket H,H\rrbracket_A=0$ if and only if $\delta_A (H^{-1})=0,$ i.e. $H^{-1}$ is a $2$-cocycle on the left-symmetric algebroid $A.$

\end{pro}

Let $(A,\cdot_A,a_A)$ be a left-symmetric algebroid, and $H\in \Sym^2(A)$. Define
\begin{equation}\label{eq:multiplication-H}
\alpha\cdot^{H^\sharp} \beta=\huaL_{H^\sharp(\alpha)}\beta-R_{H^\sharp(\beta)}\alpha-{\dM_A} (H(\alpha,\beta)), \quad\forall \alpha,\beta\in\Gamma(A^*),
\end{equation}
where $R$ is given by \eqref{eq:dualLR}, $\huaL$ and $\dM_A$ are the Lie derivation and the coboundary operator on the sub-adjacent Lie algebroid $A^c$ respectively.
\begin{pro}{\rm(\cite{lsb2})}\label{pro:morphism}
With the above notations,
for all $\alpha,\beta\in\Gamma(A^*)$, we have
\begin{equation}\label{homo}
H^\sharp(\alpha\cdot^{H^\sharp} \beta)-H^\sharp(\alpha)\cdot_A H^\sharp(\beta)=\llbracket H,H\rrbracket_A(\alpha,-,\beta).
\end{equation}
\end{pro}

\begin{thm}{\rm(\cite{lsb2})}\label{thm:LSBi-H}
With the above notations, if $\llbracket H,H\rrbracket_A=0$, then $(A^*,\cdot^{H^\sharp},a_{A^*}=a_A\circ H^\sharp)$ is a left-symmetric algebroid, and $H^\sharp$ is a left-symmetric algebroid homomorphism from $(A^*,\cdot^{H^\sharp},a_{A^*})$ to $(A,\cdot_A,a_A)$. Furthermore, $(A,A^*)$ is a left-symmetric bialgebroid.
\end{thm}

Recall that a {\bf pseudo-Hessian metric} $g$ is a
pseudo-Riemannian metric $g$ on a flat manifold $(M,\nabla)$ such
that $g$ can be locally expressed by
$g_{ij}=\frac{\partial^2\varphi}{\partial x^i\partial x^j},$ where
$\varphi\in\CWM$ and $\{x^1,\cdots,x^n\}$ is an affine coordinate
system with respect to $\nabla$. Then the pair $(\nabla,g)$ is called a
pseudo-Hessian structure on $M$. A manifold $M$ with a
pseudo-Hessian structure $(\nabla,g)$ is called a  {\bf pseudo-Hessian manifold.}
 See \cite{Geometry of Hessian structures} for more
details about pseudo-Hessian manifolds.
\begin{pro}\label{Hessian}{\rm(\cite{Geometry of Hessian structures})}
Let $(M,\nabla)$ be a flat manifold and $g$ a pseudo-Riemannian metric on $M$. Then  $g$ is a pseudo-Hessian metric if and only if for all $x,y,z\in\Gamma(TM),$ there holds $\nabla_xg(y,z)=\nabla_yg(x,z),$
where $\nabla_xg(y,z)$ is given by
\begin{eqnarray*}
\nabla_xg(y,z)=xg(y,z)-g(\nabla_x y,z)-g(y,\nabla_xz).
\end{eqnarray*}
\end{pro}
Thus we have
\begin{pro}\label{pro:hessian-cocycle}
Let $(M,\nabla)$ be a flat manifold and $g$ a pseudo-Riemannian metric on $M$. Then $(M,\nabla,g)$ is a pseudo-Hessian manifold if and only if $\delta_{T_\nabla M} g=0$, where $\delta_{T_\nabla M}$ is the coboundary operator given by \eqref{LSCA cohomology} associated to the left-symmetric algebroid $T_\nabla M $ given in Example \ref{ex:main}.
\end{pro}

Now we give the main structure studied in this paper.

\begin{defi}
Let $(A,\cdot_A,a_A)$ be a left-symmetric algebroid.
 \begin{itemize}
  \item[\rm{(i)}]  If  $H\in \Sym^2(A)$ satisfies $\llbracket H,H\rrbracket_A=0$, then $H$ is called a {\bf Koszul-Vinberg structure} on the left-symmetric algebroid $A$;
\item[\rm{(ii)}] If  $\frkB\in \Sym^2(A^*)$ satisfies $\delta_A\frkB=0$, then $\frkB$ is called a {\bf pseudo-Hessian structure} on the left-symmetric algebroid $A$.
 \end{itemize}\end{defi}

By Proposition \ref{pro:equivelent}, if a  Koszul-Vinberg structure  $H\in \Sym^2(A)$ is nondegenerate, then  $H^{-1}$ is a   pseudo-Hessian structure.


\emptycomment{
\begin{pro}
  Let $(M,\nabla)$ be a flat manifold and $\{x^1,\cdots,x^n\}$ be an affine coordinate
system with respect to $\nabla$. Then for $g\in\Sym^2(T^*M) $ and $H\in\Sym^2(TM)$,
\begin{itemize}
\item[$\rm(i)$]$g=\sum_{i,j=1}^ng_{ij}d x^i\otimes d x^j$ is a pseudo-Hessian structure on $M$ if and only if $$\frac{\partial g_{ij}}{\partial x^k}=\frac{\partial g_{kj}}{\partial x^i};$$
\item[$\rm(ii)$] $H=\sum_{i,j=1}^nh_{i,j}\frac{\partial}{\partial x^i}\otimes \frac{\partial}{\partial x^j}$ is a Koszul-Vinberg structure on $M$ if and only if $$\sum_{l=1}^n(h_{jl}\frac{\partial h_{ik}}{\partial x^l}-h_{il}\frac{\partial h_{jk}}{\partial x^l})=0.$$
\end{itemize}
\end{pro}
\begin{ex}\label{ex:affine example1}{\rm(\cite{BBo})
  Let $(\Real^n,\nabla)$ be the standard flat manifold and $\{x^1,\cdots,x^n\}$ be an affine coordinate
system with respect to $\nabla$. Let $f\in C^\infty(M)$ such that the matrix $(\frac{\partial f}{\partial x_i \partial x_j})_{r\times r}$ is nondegenerate and put
  $$H=\sum_{i,j=1}^rh_{ij}\frac{\partial}{\partial_{x_i}}\otimes \frac{\partial}{\partial_{x_j}},$$
where $(h_{ij})_{r\times r}$ is the inverse of the matrix $(\frac{\partial f}{\partial x_i \partial x_j})_{r\times r}$. Then $(\Real^n,H)$ is a Koszul-Vinberg manifold.}
\end{ex}}

Now we extend the equation \eqref{brac2} to any two elements $H_1,H_2\in\Sym^2(A)$ as follows:
\begin{eqnarray}\label{brac3}
&&\llbracket H_1,H_2\rrbracket_A(\alpha_1,\alpha_2,\alpha_3)\nonumber\\
&=&\half\Big(a_A(H^\sharp_1(\alpha_1))\langle H^\sharp_2(\alpha_2),\alpha_3\rangle+a_A(H^\sharp_2(\alpha_1))\langle H^\sharp_1(\alpha_2),\alpha_3\rangle-a_A(H^\sharp_1(\alpha_2))\langle H^\sharp_2(\alpha_1),\alpha_3\rangle\nonumber\\
&&-a_A(H^\sharp_2(\alpha_2))\langle H^\sharp_1(\alpha_1),\alpha_3\rangle+\langle \alpha_1,H^\sharp_1(\alpha_2)\cdot_A H^\sharp_2(\alpha_3)\rangle +\langle \alpha_1,H^\sharp_2(\alpha_2)\cdot_A H^\sharp_1(\alpha_3)\rangle\nonumber \\
&&-\langle\alpha_2,H^\sharp_1(\alpha_1)\cdot_A H^\sharp_2(\alpha_3)\rangle-\langle\alpha_2,H^\sharp_2(\alpha_1)\cdot_A H^\sharp_1(\alpha_3)\rangle-\langle \alpha_3,[H^\sharp_1(\alpha_1),
H^\sharp_2(\alpha_2)]_A\rangle\nonumber\\
&&-\langle \alpha_3,[H^\sharp_2(\alpha_1),
H^\sharp_1(\alpha_2)]_A\rangle\Big),\quad \forall \alpha_1,\alpha_2,\alpha_3\in\Gamma(A^*).
\end{eqnarray}
It is easy to see that $\llbracket H_1,H_2\rrbracket_A\in\Gamma(\wedge^2A\otimes A)$ and $\llbracket H_1,H_2\rrbracket_A=\llbracket H_2,H_1\rrbracket_A$.

\begin{defi}
Two Koszul-Vinberg structures $H_1$, $H_2$ are called {\bf compatible} if $\llbracket H_1,H_2\rrbracket_A=0$.
\end{defi}
It is obvious that $H_1$ and $H_2$ are compatible if and only if for all $t_1,t_2\in \K$, $t_1H_1+t_2 H_2$ are Koszul-Vinberg structures.

\section{Koszul-Vinberg-Nijenhuis structures on left-symmetric algebroids}\label{sec:KVN}

Let $N$ be a Nijenhuis operator on a left-symmetric algebroid $({A},\cdot_{A},a_A)$.
For all $x\in\Gamma(A)$, we denote by
${L}^N_x:\Gamma(A)\longrightarrow\Gamma(A)$  and $R^N_x:\Gamma(A)\longrightarrow\Gamma(A)$ the left multiplication and the right multiplication  for the left-symmetric algebroid $({A},\cdot_N,a_N)$ respectively, i.e.,
\begin{equation*}\label{eq:leftmu2}{L}^N_xy=x\cdot_N y,\quad{R}^N_xy=y\cdot_N x.\end{equation*}
By direct calculations, for all $\alpha\in\Gamma(A^*)$, we have
\begin{eqnarray}
 {L}^N_x\alpha&=&L_{N(x)}\alpha+N^*(L_{x}\alpha)-L_{x}N^*(\alpha),\label{eq:pre-Lie leftmulti1}\\
{R}^N_x\alpha&=&R_{N(x)}\alpha+N^*(R_{x}\alpha)-R_{x}N^*(\alpha).\label{eq:pre-Lie leftmulti2}
\end{eqnarray}
For all $x\in\Gamma(A)$, the Lie derivation for the Lie algebroid $({A},[-,-]_N,a_N)$ is denoted by $\huaL^N_x$, i.e.,
\begin{equation*}
\langle\huaL^N_x\alpha,y\rangle=a_N(x)\langle\alpha,y\rangle-\langle \alpha,[x,y]_N\rangle.
\end{equation*}
By direct calculations, we have
\begin{eqnarray}
 \huaL^N_x\alpha&=&\huaL_{N(x)}\alpha+N^*(\huaL_{x}\alpha)-\huaL_{x}N^*(\alpha).\label{eq:pre-Lie leftmulti3}
\end{eqnarray}
Let $H\in\Sym^2(A)$ be a Koszul-Vinberg  structure. We define the multiplication $\cdot_{N^*}^{H^\sharp}:\Gamma({A}^*)\otimes\Gamma({A}^*)\longrightarrow\Gamma({A}^*)$ to be the deformed multiplication of $\cdot^{H^\sharp}$ by $N^*$. More precisely,
\begin{eqnarray}\label{eq:KVN-operation1}
\alpha\cdot_{N^*}^{H^\sharp}\beta=N^*(\alpha)\cdot^{H^\sharp} \beta+\alpha\cdot^{H^\sharp}N^*(\beta)-N^*(\alpha\cdot^{H^\sharp}\beta),\quad\forall~\alpha,\beta\in\Gamma({A}^*).
\end{eqnarray}
Define the multiplication $\star^{H^\sharp}:\Gamma({A}^*)\otimes\Gamma({A}^*)\longrightarrow\Gamma({A}^*)$ similar as \eqref{eq:multiplication-H} by
\begin{eqnarray}
  \alpha\star^{H^\sharp }\beta&=&\huaL^N_{H^\sharp(\alpha)}\beta-R^N_{{H^\sharp(\beta)}}\alpha-{\dM^N_A} (H(\alpha,\beta)),
\end{eqnarray}
where $\huaL^N$ is given by \eqref{eq:pre-Lie leftmulti3}, $R^N$ is given by \eqref{eq:pre-Lie leftmulti2}, and ${\dM^N_A}$ is the coboundary operator associated to the Lie algebroid $({A},[-,-]_N,a_N)$.
\begin{lem}\label{lem:threeformular1}
With the above notations, if $N\circ H^\sharp=H^\sharp\circ N^*$, we have
\begin{equation}
  \alpha\star^{H^\sharp}\beta+\alpha\cdot_{N^*}^{H^\sharp}\beta=2\alpha\cdot^{N\circ H^\sharp}\beta,\quad \alpha, \beta\in\Gamma(A^*).
\end{equation}
\end{lem}
\pf Since $N\circ H^\sharp=H^\sharp\circ N^*$, by direct calculations, we have
\begin{eqnarray*}
   \alpha\star^{H^\sharp}\beta+\alpha\cdot_{N^*}^{H^\sharp}\beta&=&\huaL_{N\circ H^\sharp(\alpha)}\beta+N^*(\huaL_{H^\sharp(\alpha)}\beta)-\huaL_{H^\sharp(\alpha)}N^*(\beta)-R_{N\circ H^\sharp(\beta)}\alpha\\
   &&-N^*(R_{H^\sharp(\beta)}\alpha)+R_{H^\sharp(\beta)}N^*(\alpha)-N^*\dM_AH(\alpha,\beta)+\huaL_{H^\sharp\circ N^*(\alpha)}\beta\\
   &&-R_{H^\sharp(\beta)}N^*(\alpha)-{\dM_A} H(N^*(\alpha),\beta)+\huaL_{H^\sharp(\alpha)}N^*(\beta)-R_{H^\sharp\circ N^*(\beta)}\alpha\\
   &&-{\dM_A} H(\alpha,N^*(\beta))-N^*\big(\huaL_{H^\sharp(\alpha)}\beta-R_{H^\sharp(\beta)}\alpha-\dM_AH(\alpha,\beta)\big)\\
   &=&\huaL_{N\circ H^\sharp(\alpha)}\beta+\huaL_{H^\sharp\circ N^*(\alpha)}\beta-R_{N\circ H^\sharp(\beta)}\alpha-R_{H^\sharp\circ N^*(\beta)}\alpha\\
   &&-{\dM_A} H(N^*(\alpha),\beta)-{\dM_A} H(\alpha,N^*(\beta))\\
   &=&2\alpha\cdot^{N\circ H^\sharp}\beta.
\end{eqnarray*}
We finish the proof.\qed\vspace{3mm}

We denote the commutators of the operations $\star^{H^\sharp}$ and $\cdot_{N^*}^{H^\sharp}$ by $\{-,-\}^{H^\sharp}$ and $[-,-]_{N^*}^{H^\sharp}$ respectively:
\begin{eqnarray}
 \{\alpha,\beta\}^{H^\sharp}&=&\alpha\star^{H^\sharp}\beta-\beta\star^{H^\sharp}\alpha=L^N_{H^\sharp(\alpha)}\beta-L^N_{H^\sharp(\beta)}\alpha,\label{eq:sub-ad Lie 1}\\
 {[\alpha,\beta]_{N^*}^{H^\sharp}}&=&\alpha\cdot_{N^*}^{H^\sharp}\beta-\beta\cdot_{N^*}^{H^\sharp}\alpha=[N^*(\alpha),\beta]^{{H^\sharp}}+[\alpha,N^*(\beta)]^{H^\sharp}-N^*([\alpha,\beta]^{H^\sharp}),\label{eq:sub-ad Lie 2}
\end{eqnarray}
where $[-,-]^{H^\sharp}$ is the commutator of the operation $\cdot^{H^\sharp}$, i.e.
\begin{equation}\label{eq:commutator-H}
  [\alpha,\beta]^{H^\sharp}=L_{H^\sharp(\alpha)}\beta-L_{H^\sharp(\beta)}\alpha,\quad\forall~\alpha,\beta\in\Gamma(A^*).
\end{equation}

In the following, we introduce the definition of a Koszul-Vinberg-Nijenhuis structure on a left-symmetric algebroid.
\begin{defi}
Let $H$ be a Koszul-Vinberg structure and $N:A\longrightarrow A$ a Nijenhuis operator on a left-symmetric algebroid $({A},\cdot_{A},a_A)$ . A pair $(H,N)$ is called a {\bf Koszul-Vinberg-Nijenhuis structure} (or {\bf $\KVN$-structure}) on the left-symmetric algebroid ${A}$ if for any $\alpha,\beta\in\Gamma({A}^*)$,
\begin{eqnarray}
 \label{eq:BNS1} N\circ H^\sharp&=&H^\sharp\circ N^*,\\
 \label{eq:BNS2} \alpha\star^{H^\sharp}\beta&=&\alpha\cdot_{N^*}^{H^\sharp}\beta.
\end{eqnarray}
\end{defi}
By Lemma \ref{lem:threeformular1}, if $(H,N)$ is a $\KVN$-structure, we have
\begin{equation}\label{eq:BSN4}
\alpha\star^{H^\sharp}\beta=\alpha\cdot_{N^*}^{H^\sharp}\beta=\alpha\cdot^{N\circ H^\sharp}\beta.
\end{equation}

\begin{thm}\label{thm:properties of GHN}
 Let $(H,N)$ be a $\KVN$-structure on a left-symmetric algebroid $({A},\cdot_{A},a_A)$. Then
 \begin{itemize}
\item[$\rm(i)$]$N^*$ is a Nijenhuis operator on the left-symmetric algebroid $({A}^*,\cdot^{H^\sharp},a_{A^*}=a_A\circ H^\sharp)$. Consequently,  $({A}^*,\cdot^{H^\sharp}_ {N^*},a_{NH}=a_A\circ N\circ H^\sharp)$ is also a left-symmetric algebroid and $N^*$ is a left-symmetric algebroid morphism from $({A^*},\cdot^{H^\sharp}_ {N^*},a_{NH})$ to  $({A^*},\cdot^{H^\sharp},a_{A^*})$;
\item[$\rm(ii)$]$H^\sharp$ is a left-symmetric algebroid morphism from $({A}^*,\cdot^{H^\sharp}_ {N^*},a_{NH})$ to $({A},\cdot_N,a_N)$;
\item[$\rm(iii)$]$N\circ H^\sharp$ is a left-symmetric algebroid morphism from $({A^*},\cdot^{H^\sharp}_ {N^*},a_{NH})$ to $({A},\cdot_{A},a_A)$;
\item[$\rm(iv)$] $H_N\in\Sym^2(A)$, which is defined by $(H_N)^\sharp=N\circ H^\sharp$, satisfies $\llbracket   H_N, H_N\rrbracket_A=0$ and $H$ is also a Koszul-Vinberg structure for the deformed left-symmetric algebroid $({A},\cdot_N,a_N)$;
\item[$\rm(v)$]$\llbracket H,H_N\rrbracket_A=0$.
\end{itemize}

\end{thm}
\pf (i) By \eqref{eq:BNS2}, we have
\begin{eqnarray}\label{eq:operations1}
R_{H^\sharp(\beta)}N^*(\alpha)=\huaL_{H^\sharp(\alpha)}N^*(\beta)-N^*\huaL_{H^\sharp(\alpha)}\beta+N^*R_{H^\sharp(\beta)}\alpha+N^*\dM_AH(\alpha,\beta)-\dM_AH(N^*(\alpha),\beta).
\end{eqnarray}
Replacing $\alpha$ by $N^*(\alpha)$ in \eqref{eq:operations1}, and by \eqref{eq:BNS1}, we have
\begin{eqnarray}\label{eq:operations2}
R_{H^\sharp(\beta)}(N^*)^2(\alpha)&=&\huaL_{NH^\sharp(\alpha)}N^*(\beta)-N^*\huaL_{NH^\sharp(\alpha)}\beta+N^*R_{H^\sharp(\beta)}N^*(\alpha)\nonumber\\
&&+N^*\dM_AH(N^*(\alpha),\beta)-\dM_AH((N^*)^2(\alpha),\beta).
\end{eqnarray}
Let $N^*$ act both sides of \eqref{eq:operations1}, we have
\begin{eqnarray}\label{eq:operations3}
N^*R_{H^\sharp(\beta)}N^*(\alpha)&=&N^*\huaL_{H^\sharp(\alpha)}N^*(\beta)-(N^*)^2\huaL_{H^\sharp(\alpha)}\beta+(N^*)^2R_{H^\sharp(\beta)}\alpha\nonumber\\
&&+(N^*)^2\dM_AH(\alpha,\beta)-N^*\dM_AH(N^*(\alpha),\beta).
\end{eqnarray}
Since $N$ is a Nijenhuis operator on the left-symmetric algebroid $A$, we have
\begin{eqnarray}\label{eq:operations4}
R_{H^\sharp(\beta)}(N^*)^2(\alpha)-N^*R_{H^\sharp(\beta)}N^*(\alpha)-R_{H^\sharp N^*(\beta)}\alpha+N^*R_{H^\sharp N^*(\beta)}\alpha=0.
\end{eqnarray}
By \eqref{eq:BNS1} and \eqref{eq:operations2}-\eqref{eq:operations4}, we have
\begin{eqnarray*}
&&N^*(\alpha)\cdot^{H^\sharp} N^*(\beta)-N^*(\alpha\cdot^{H^\sharp}_{N^*}\beta)\\
&=&\huaL_{NH^\sharp(\alpha)}N^*(\beta)-N^*\huaL_{NH^\sharp(\alpha)}\beta+N^*R_{H^\sharp(\beta)}N^*(\alpha)+N^*\dM_AH(N^*(\alpha),\beta)\\
&&-\dM_AH((N^*)^2(\alpha),\beta)+(N^*)^2\huaL_{H(\alpha)}\beta-N^*\huaL_{H(\alpha)}N^*(\beta)-(N^*)^2R_{H(\beta)}\alpha\\
&&-(N^*)^2\dM_AH(\alpha,\beta)+N^*\dM_AH(N^*(\alpha),\beta)-R_{NH^\sharp(\beta)}N^*(\alpha)+N^*R_{NH^\sharp(\beta)}\alpha\\
&=&R_{H^\sharp(\beta)}(N^*)^2(\alpha)-N^*R_{H^\sharp(\beta)}N^*(\alpha)-R_{NH^\sharp(\beta)}N^*(\alpha)+N^*R_{NH^\sharp(\beta)}\alpha=0.
\end{eqnarray*}
Thus $N^*$ is a Nijenhuis operator on the left-symmetric algebroid $({A}^*,\cdot^{H^\sharp},a_{A^*}=a_A\circ H^\sharp)$. The rest follows from the properties of Nijenhuis operators.

(ii) By \eqref{eq:BNS1} and Theorem \ref{thm:LSBi-H}, we have
\begin{eqnarray*}
 H^\sharp( \alpha\cdot_{N^*}^{H^\sharp}\beta)&=&H^\sharp\big(N^*(\alpha)\cdot^{H^\sharp} \beta+\alpha\cdot^{H^\sharp}N^*(\beta)-N^*(\alpha\cdot^{H^\sharp}\beta)\big)\\
  &=&H^\sharp N^*(\alpha)\cdot_A H^\sharp(\beta)+H^\sharp(\alpha)\cdot_AH^\sharp N^*(\beta)-H^\sharp N^*(\alpha\cdot^{H^\sharp}\beta)\\
  &=&NH^\sharp(\alpha)\cdot_A H^\sharp(\beta)+H^\sharp(\alpha)\cdot_ANH^\sharp (\beta)-N(H^\sharp(\alpha)\cdot_AH^\sharp(\beta))\\
  &=&H^\sharp(\alpha)\cdot_N H^\sharp(\beta).
\end{eqnarray*}
Also, $a_{NH}=a_N\circ H^\sharp$, thus $H^\sharp$ is a left-symmetric algebroid morphism from $({A},\cdot^{H}_ {N^*},a_{NH})$ to $({A},\cdot_N,a_N)$.

(iii) By (ii) and $N$ is a Nijenhuis operator, we have
\begin{eqnarray*}
 NH^\sharp(\alpha\cdot_{N^*}^{H^\sharp}\beta)=N(H^\sharp(\alpha)\cdot_N H^\sharp(\beta))=NH^\sharp(\alpha)\cdot_ANH^\sharp(\beta).
\end{eqnarray*}
Also, $a_{NH}=a_A\circ N\circ H^\sharp$, thus $N\circ H^\sharp$ is a left-symmetric algebroid morphism from $({A^*},\cdot^{H^\sharp}_ {N^*},a_{NH})$ to $({A},\cdot_{A},a_A)$.

(iv) By (iii) and \eqref{eq:BSN4}, we have
\begin{eqnarray*}
   NH^\sharp(\alpha\cdot^{N\circ H^\sharp}\beta)= NH^\sharp(\alpha\cdot_{N^*}^{H^\sharp}\beta)=NH^\sharp(\alpha)\cdot_ANH^\sharp(\beta).
\end{eqnarray*}
By Proposition \ref{pro:morphism}, we have $\llbracket H_N,H_N\rrbracket_A=0$.

By (ii) and \eqref{eq:BSN4}, we have
\begin{eqnarray*}
    H^\sharp(\alpha\star^{H^\sharp}\beta)= H^\sharp( \alpha\cdot_{N^*}^{H^\sharp}\beta)=H^\sharp(\alpha)\cdot_NH^\sharp(\beta).
\end{eqnarray*}
Also, by Proposition \ref{pro:morphism}, $H$ is a Koszul-Vinberg structure for the deformed left-symmetric algebroid $({A},\cdot_N,a_N)$.

(v) First, we notice that
\begin{eqnarray*}
  \alpha\cdot^{H^\sharp+N\circ H^\sharp}\beta=\alpha\cdot^{H^\sharp}\beta+\alpha\cdot^{N\circ H^\sharp}\beta=\alpha\cdot^{H^\sharp}\beta+\alpha\cdot_{H^\sharp}^{N^*}\beta.
\end{eqnarray*}
Then we have
\begin{eqnarray*}
&&(H+H_N)^\sharp(\alpha\cdot^{(H+H_N)^\sharp}\beta)\\
&=&(H^\sharp+N\circ H^\sharp) (\alpha\cdot^{H^\sharp+N\circ H^\sharp}\beta)\\
 &=&H^\sharp(\alpha\cdot^{H^\sharp}\beta)+H^\sharp(\alpha\cdot^{H^\sharp}_{N^*}\beta)+NH^\sharp(\alpha\cdot^{H^\sharp}\beta)+NH^\sharp(\alpha\cdot^{H^\sharp}_{N^*}\beta)\\
 &=&H^\sharp(\alpha\cdot^{H^\sharp}\beta)+H^\sharp\big(N^*(\alpha)\cdot^{H^\sharp} \beta+\alpha\cdot^{H^\sharp}N^*(\beta)-N^*(\alpha\cdot^{H^\sharp}\beta)\big)\\
 &&+NH^\sharp(\alpha\cdot^{H^\sharp}\beta)+NH^\sharp(\alpha\cdot^{H^\sharp}_{N^*}\beta)\\
 &=&H^\sharp(\alpha\cdot^{H^\sharp}\beta)+H^\sharp(N^*(\alpha)\cdot^{H^\sharp} \beta)+H^\sharp(\alpha\cdot^{H^\sharp}N^*(\beta))+NH^\sharp(\alpha\cdot^{H^\sharp}_{N^*}\beta)\\
 &=&H^\sharp(\alpha)\cdot_AH^\sharp(\beta)+NH^\sharp(\alpha)\cdot_AH^\sharp(\beta)+H^\sharp(\alpha)\cdot_A NH^\sharp(\beta)+NH^\sharp(\alpha)\cdot_A NH^\sharp(\beta)\\
 &=& (H^\sharp+N\circ H^\sharp)(\alpha)\cdot_A(H^\sharp+N\circ H^\sharp)(\beta)\\
 &=&(H+H_N)^\sharp(\alpha)\cdot_A(H+H_N)^\sharp(\beta).
\end{eqnarray*}
Thus, by Proposition \ref{pro:morphism}, we have $\llbracket H+H_N,H+H_N\rrbracket_A=0$ and then $\llbracket H,H_N\rrbracket_A=0$.\qed

\begin{cor}
  Let $(H,N)$ be a $\KVN$-structure on a left-symmetric algebroid $({A},\cdot_{A},a_A)$. Then
 \begin{itemize}
\item[$\rm(i)$]$N^*$ is a Nijenhuis operator on the Lie algebroid $({A}^*,[-,-]^{H^\sharp},a_{A^*}=a_A\circ H^\sharp)$. Consequently,  $({A^*},[-,-]^{H^\sharp}_ {N^*},a_{NH}=a_A\circ N\circ H^\sharp)$ is  a Lie algebroid and $N^*$ is a Lie algebroid morphism from $({A^*},[-,-]^{H^\sharp}_ {N^*},a_{NH})$ to   $({A^*},[-,-]^{H^\sharp},a_{A^*})$;
\item[$\rm(ii)$]$H^\sharp$ is a Lie algebroid morphism from $({A^*},[-,-]^{H^\sharp}_ {N^*},a_{NH})$ to $({A},[-,-]_N,a_N)$;
\item[$\rm(iii)$]$N\circ H^\sharp$ is a Lie algebroid morphism from $({A^*},[-,-]^{H^\sharp}_ {N^*},a_{NH})$ to $({A},[-,-]_{A},a_A)$.
\end{itemize}
\end{cor}
\section{Hierarchy of Koszul-Vinberg structures}\label{sec:comgRB}
Given a $\KVN$-structure $(H,N)$, by Theorem \ref{thm:properties of GHN}, $H$ and $H_N$ are Koszul-Vinberg structures and they are compatible. In the following, we construct a hierarchy of Koszul-Vinberg structures from a $\KVN$-structure. We first give two useful lemmas. Define $H_{N^k}\in\Sym^2(A)$ by $(H_{N^k})^\sharp=N^k\circ H^\sharp$ for all $k\in \mathbb N.$

\begin{lem}
  Let $(H,N)$ be a $\KVN$-structure on a left-symmetric algebroid $(A,\cdot_A,a_A)$. Then for all $k,i\in\Nat$, we have
\begin{eqnarray}
   \label{eq:HSN3}(H_{N^k})^\sharp(\alpha\cdot^{H^\sharp}_{{({N^*})}^{k+i}}\beta)= (H_{N^k})^\sharp(\alpha)\cdot_{N^i}(H_{N^k})^\sharp(\beta),\quad \forall \alpha,\beta\in\Gamma(A^*).
\end{eqnarray}

\end{lem}
\pf Since $H$ is a Koszul-Vinberg structure and $H^\sharp\circ{N^*}=N\circ H^\sharp$, we have
\begin{eqnarray}
\nonumber H^\sharp(\alpha\cdot^{H^\sharp}_{{({N^*})}^i}\beta)&=&H^\sharp\big({({N^*})}^i(\alpha)\cdot^{H^\sharp} \beta+\alpha\cdot^{H^\sharp} {({N^*})}^i(\beta)-{({N^*})}^{i}(\alpha\cdot^{H^\sharp} \beta)\big)\\
\nonumber&=&N^i(H^\sharp (\alpha))\cdot_A H^\sharp(\beta)+H^\sharp(\alpha)\cdot_A N^i(H^\sharp(\beta))-N^i(H^\sharp(\alpha)\cdot_A H^\sharp(\beta))\\
\label{eq:HSN2}&=& H^\sharp(\alpha)\cdot_{N^i} H^\sharp(\beta).
\end{eqnarray}
Since $N^*$ is a Nijenhuis operator on the left-symmetric algebroid $({A}^*,\cdot^{H^\sharp},a_{A^*}=a_A\circ H^\sharp)$, by Lemma \ref{lem:Niejproperty}, we have
\begin{equation}\label{eq:t1}
  {({N^*})}^k(\alpha\cdot^{H^\sharp}_{{({N^*})}^{k+i}}\beta)={({N^*})}^k(\alpha)\cdot^{H^\sharp}_{{({N^*})}^i} {({N^*})}^k(\beta).
\end{equation}
By \eqref{eq:HSN2} and \eqref{eq:t1}, we have
\begin{eqnarray*}
(H_{N^k})^\sharp(\alpha\cdot^{H^\sharp}_{{({N^*})}^{k+i}}\beta)&=&H^\sharp\circ {({N^*})}^k(\alpha\cdot^{H^\sharp}_{{({N^*})}^{k+i}}\beta)\\
 &=&H^\sharp\big({({N^*})}^k(\alpha)\cdot^{H^\sharp}_{{({N^*})}^i}{({N^*})}^k(\beta)\big)\\
 &=& H^\sharp({({N^*})}^k(\alpha))\cdot_{N^i}H^\sharp({({N^*})}^k(\beta)),\\
 &=&(H_{N^k})^\sharp(\alpha)\cdot_{N^i}(H_{N^k})^\sharp(\beta).
\end{eqnarray*}
The proof is finished. \qed

\begin{lem}
 Let $(H,N)$ be a $\KVN$-structure on a left-symmetric algebroid $(A,\cdot_A,a_A)$. Then for all $k,i\in\Nat$ such that $i\leq k$,
\begin{equation}\label{eq:HSN4}
 \alpha\cdot^{(H_{N^k})^\sharp}\beta= \alpha\cdot^{H^\sharp}_{{({N^*})}^k}\beta={({N^*})}^{k-i}(\alpha\cdot^{(H_{N^i})^\sharp}\beta).
\end{equation}
\end{lem}
\pf
 Since ${N}$ and $H$ are compatible, the relation $ \alpha\cdot^{(H_{N^k})^\sharp}\beta= \alpha\cdot^{H^\sharp}_{{({N^*})}^k}\beta$ is valid for $k=0$ and $k=1$. Let us assume that this relation is valid for all integers less than or equal to $k\geq1$. By direct calculation, we have
\begin{eqnarray*}
  \alpha\cdot^{(H_{N^{k+1}})^\sharp}\beta-\alpha\cdot^{(H_{N^{k}})^\sharp}_{N^*} \beta-\big({N^*}(\alpha\cdot^{(H_{N^k})^\sharp}\beta)-{N^*}(\alpha)\cdot^{(H_{N^{k-1}})^\sharp}{N^*}(\beta)\big)=0.
\end{eqnarray*}
By \eqref{eq:t1}, we have
\begin{eqnarray*}
{N^*}(\alpha\cdot^{(H_{N^k})^\sharp}\beta)-{N^*}(\alpha)\cdot^{(H_{N^{k-1}})^\sharp}{N^*}(\beta)={N^*}(\alpha\cdot^{H^\sharp}_{{({N^*})}^k}\beta)-{N^*}(\alpha)\cdot^{H^\sharp}_{{({N^*})}^{k-1}}{N^*}(\beta)=0.
\end{eqnarray*}
Thus
$$\alpha\cdot^{(H_{N^{k+1}})^\sharp}\beta=\alpha\cdot_{{N^*}}^{(H_{N^k})^\sharp}\beta={N^*}(\alpha\cdot^{(H_{N^{k}})^\sharp}\beta)={N^*}(\alpha\cdot^{H^\sharp}_{{({N^*})}^k}\beta)=\alpha\cdot^{H^\sharp}_{{({N^*})}^{k+1}}\beta.$$
 The relation $ \alpha\cdot^{(H_{N^k})^\sharp}\beta= \alpha\cdot^{H^\sharp}_{{({N^*})}^k}\beta$ is proved by induction on $k$. The second equality in \eqref{eq:HSN4} follows from the relations $ \alpha\cdot^{(H_{N^k})^\sharp}\beta= \alpha\cdot^{H^\sharp}_{{({N^*})}^k}\beta$ and $\alpha\cdot_{{N^*}}^{(H_{N^k})^\sharp}\beta={N^*}(\alpha\cdot^{(H_{N^{k}})^\sharp}\beta)$. We finish the proof.\qed\vspace{3mm}

 With above preparations, we give the main result in this section.

\begin{thm}\label{thm:hierarchy}
 Let $(H,N)$ be a $\KVN$-structure on a left-symmetric algebroid $(A,\cdot_A,a_A)$. Then for all $k\in\Nat$, $H_{N^{k}}$ are Koszul-Vinberg structures. Furthermore, for all $k,l\in\Nat$, $H_{N^{k}}$ and $H_{N^{l}}$ are compatible.
\end{thm}
\pf By  \eqref{eq:HSN3} and \eqref{eq:HSN4} with $k=0$, we have
\begin{eqnarray*}
 (H_{N^k})^\sharp(\alpha\cdot^{(H_{N^{k}})^\sharp}\beta)= (H_{N^k})^\sharp(\alpha)\cdot_A (H_{N^k})^\sharp(\beta),
\end{eqnarray*}
which implies that $H_{N^{k}}$ is a Koszul-Vinberg structure. For the second conclusion, we need  to prove that $H_{N^{k}}+H_{N^{k+i}}$ are Koszul-Vinberg structures for all $k,i\in\Nat$.
By  \eqref{eq:HSN4}, we have
$$\alpha\cdot^{(H_{N^k})^\sharp+(H_{N^{k+i}})^\sharp}\beta=\alpha\cdot^{(H_{N^k})^\sharp}\beta+\alpha\cdot^{(H_{N^{k+i}})^\sharp}\beta=\alpha\cdot^{(H_{N^k})^\sharp}\beta+\alpha\cdot^{(H_{N^k})^\sharp}_{{({N^*})}^i}\beta.$$
Thus, we have
\begin{eqnarray*}
  &&((H_{N^k})^\sharp+(H_{N^{k+i}})^\sharp)(\alpha\cdot^{(H_{N^k})^\sharp+(H_{N^{k+i}})^\sharp}\beta)\\
  &=&(H_{N^k})^\sharp(\alpha\cdot^{(H_{N^k})^\sharp}\beta)+(H_{N^k})^\sharp(\alpha\cdot^{(H_{N^k})^\sharp}_{{({N^*})}^i}\beta)+(H_{N^{k+i}})^\sharp(\alpha\cdot^{(H_{N^k})^\sharp}\beta)+(H_{N^{k+i}})^\sharp(\alpha\cdot^{(H_{N^k})^\sharp}_{{({N^*})}^i}\beta)\\
  &=&(H_{N^k})^\sharp(\alpha\cdot^{(H_{N^k})^\sharp}\beta)+(H_{N^{k+i}})^\sharp(\alpha\cdot^{(H_{N^k})^\sharp}\beta)+(H_{N^{k+i}})^\sharp(\alpha\cdot^{(H_{N^k})^\sharp}_{{({N^*})}^i}\beta)\\
  &&+(H_{N^k})^\sharp\big({({N^*})}^i(\alpha)\cdot^{(H_{N^k})^\sharp}\beta+\alpha\cdot^{(H_{N^k})^\sharp}{({N^*})}^i(\beta)-{({N^*})}^i(\alpha\cdot^{(H_{N^k})^\sharp}\beta)\big)\\
  &=&(H_{N^k})^\sharp(\alpha\cdot^{(H_{N^k})^\sharp}\beta)+(H_{N^{k+i}})^\sharp(\alpha\cdot^{(H_{N^k})^\sharp}_{{({N^*})}^i}\beta)+(H_{N^k})^\sharp\big({({N^*})}^i(\alpha)\cdot^{(H_{N^k})^\sharp}\beta+\alpha\cdot^{(H_{N^k})^\sharp}{({N^*})}^i(\beta)\big)\\
  &=&(H_{N^k})^\sharp(\alpha)\cdot_A (H_{N^k})^\sharp(\beta)+(H_{N^{k+i}})^\sharp(\alpha)\cdot_A (H_{N^{k+i}})^\sharp(\beta)+(H_{N^{k+i}})^\sharp(\alpha)\cdot_A (H_{N^k})^\sharp(\beta)\\
  &&+(H_{N^k})^\sharp(\alpha)\cdot_A (H_{N^{k+i}})^\sharp(\beta)\\
  &=&((H_{N^k})^\sharp+(H_{N^{k+i}})^\sharp)(\alpha)\cdot_A ((H_{N^k})^\sharp+(H_{N^{k+i}})^\sharp)(\beta).
\end{eqnarray*}
Thus $H_{N^{k}}+H_{N^{k+i}}$ is a $\KVN$-structure. We finish the proof.\qed\vspace{3mm}

Compatible Koszul-Vinberg structures can give rise to $\KVN$-structures.
\begin{thm}\label{thm:Comp-GHN}
  Let $H$ and $H_1$ be two Koszul-Vinberg structures on a left-symmetric algebroid $(A,\cdot_A,a_A)$. Suppose that $H$ is nondegenerate. If $H$ and $H_1$ are compatible, then
  \begin{itemize}
  \item[$\rm(i)$] $N=H^\sharp_1\circ (H^\sharp)^{-1}$ is a Nijenhuis operator on the left-symmetric algebroid $(A,\cdot_A,a_A)$;
\item[$\rm(ii)$]$(H,N)$ is a $\KVN$-structure;
    \item[$\rm(iii)$]$(H_1,N)$ is a $\KVN$-structure.
\end{itemize}
\end{thm}

\pf (i) For all $x,y\in \Gamma(A)$, there exists $\alpha,\beta\in \Gamma(A^*)$ such that
$H^\sharp(\alpha)=x, H^\sharp(\beta)=y$. Hence $N=H^\sharp_1\circ (H^\sharp)^{-1}$ is a Nijenhuis operator
if and only if the following equation holds:
$$NH^\sharp(\alpha)\cdot_A NH^\sharp(\beta )=N\big(NH^\sharp(\alpha)\cdot_A H^\sharp(\beta )+H^\sharp(\alpha)\cdot_A NH^\sharp(\beta )\big)-N^2(H^\sharp(\alpha)\cdot_A H^\sharp(\beta )).$$
Since $H_1$ is a Koszul-Vinberg structure with $H^\sharp_1=N\circ H^\sharp$, the left hand side of the above equation is
$$NH^\sharp(\huaL_{NH^\sharp(\alpha)}\beta-R_{NH^\sharp(\beta)}\alpha-\dM_A\langle NH^\sharp(\alpha),\beta\rangle).$$
Since $H$ and $H_1$ are compatible Koszul-Vinberg structures, we have
\begin{eqnarray*}
 &&NH^\sharp(\alpha)\cdot_A H^\sharp(\beta )+H^\sharp(\alpha)\cdot_A NH^\sharp(\beta )\\
 &=&NH^\sharp\big(\huaL_{H^\sharp(\alpha)}\beta-R_{H^\sharp(\beta)}\alpha-\dM_A\langle H^\sharp(\alpha),\beta\rangle\big)+H^\sharp\big(\huaL_{NH^\sharp(\alpha)}\beta-R_{NH^\sharp(\beta)}\alpha-\dM_A\langle NH^\sharp(\alpha),\beta\rangle\big)\\
&=&N(H^\sharp(\alpha)\cdot_A H^\sharp(\beta ))+H^\sharp\big(\huaL_{NH^\sharp(\alpha)}\beta-R_{NH^\sharp(\beta)}\alpha-\dM_A\langle NH^\sharp(\alpha),\beta\rangle\big).
\end{eqnarray*}
Let $N$ act on both sides, we get the conclusion.

(ii) On the one hand, since $H$ and $H_1$ are compatible Koszul-Vinberg structures with $H^\sharp_1=N\circ H^\sharp$,
\begin{eqnarray*}
 NH^\sharp(\alpha)\cdot_A H^\sharp(\beta )+H^\sharp(\alpha)\cdot_A NH^\sharp(\beta )
 &=&NH^\sharp\big(\huaL_{H^\sharp(\alpha)}\beta-R_{H^\sharp(\beta)}\alpha-\dM_A\langle H^\sharp(\alpha),\beta\rangle\big)\\
 &&+H^\sharp\big(\huaL_{NH^\sharp(\alpha)}\beta-R_{NH^\sharp(\beta)}\alpha-\dM_A\langle NH^\sharp(\alpha),\beta\rangle\big).
 \end{eqnarray*}
On the other hand, since $H_1$ is a Koszul-Vinberg structure, we have
\begin{eqnarray*}
 &&NH^\sharp(\alpha)\cdot_A H^\sharp(\beta )+H^\sharp(\alpha)\cdot_A NH^\sharp(\beta )\\
 &=& H^\sharp N^*(\alpha)\cdot_A H^\sharp(\beta)+H^\sharp(\alpha)\cdot_AH^\sharp N^*(\beta )\\
 &=&H^\sharp\big(\huaL_{H^\sharp N^*(\alpha)}\beta-R_{H^\sharp(\beta)}N^*(\alpha)-\dM_A\langle H^\sharp N^*(\alpha),\beta\rangle+\huaL_{H^\sharp(\alpha)}N^*(\beta)-R_{H^\sharp N^*(\beta)}\alpha\\
 &&-\dM_A\langle H^\sharp (\alpha),N^*(\beta)\rangle \big)
\end{eqnarray*}
Compare the two equalities above with $H$ nondegenerate, we have
\begin{eqnarray}\label{eq:Compatible-HN1}
N^*\big(\huaL_{H^\sharp(\alpha)}\beta-R_{H^\sharp(\beta)}\alpha-\dM_A\langle H^\sharp(\alpha),\beta\rangle\big)
=\huaL_{H^\sharp(\alpha)}N^*(\beta)-R_{H^\sharp(\beta)}N^*(\alpha)-\dM_A\langle H^\sharp (\alpha),N^*(\beta)\rangle.
\end{eqnarray}
By direct calculations, we have
\begin{eqnarray*}
\alpha\star^{H^\sharp}\beta-\alpha\cdot_{N^*}^{H^\sharp}\beta&=&N^*\huaL_{H^\sharp(\alpha)}\beta-\huaL_{H^\sharp(\alpha)}N^*(\beta)-N^*R_{H^\sharp(\beta)}\alpha\\
&&+R_{H^\sharp(\beta)}N^*(\alpha)-N^*\dM_AH(\alpha,\beta)+\dM_AH(N^*(\alpha),\beta),
\end{eqnarray*}
which is just the equality \eqref{eq:Compatible-HN1}. Thus, $(H,N)$ is a $\KVN$-structure.

(iii) By direct calculations, we have
\begin{eqnarray*}
&&\alpha\star^{H_1^\sharp}\beta-\alpha\cdot_{N^*}^{H_1^\sharp}\beta\\
&=&N^*\huaL_{H_1^\sharp(\alpha)}\beta-\huaL_{H_1^\sharp(\alpha)}N^*(\beta)-N^*R_{H_1^\sharp(\beta)}\alpha+R_{H_1^\sharp(\beta)}N^*(\alpha)-N^*\dM_AH_1(\alpha,\beta)\\
&&+\dM_AH_1(N^*(\alpha),\beta)\\
&=&N^*\huaL_{NH^\sharp(\alpha)}\beta-\huaL_{NH^\sharp(\alpha)}N^*(\beta)-N^*R_{NH^\sharp(\beta)}\alpha+R_{NH^\sharp(\beta)}N^*(\alpha)-N^*\dM_A\langle NH^\sharp(\alpha),\beta\rangle\\
&&+\dM_A\langle N^2H^\sharp(\alpha),\beta\rangle\\
&=&N^*(\alpha\cdot^{N\circ H^\sharp} \beta)-N^*(\alpha)\cdot^{H^\sharp} N^*(\beta)=0,
\end{eqnarray*}
which implies that $(H_1,N)$ is a $\KVN$-structure.\qed\vspace{3mm}

The above theorem admits a converse.

\begin{thm}\label{thm:Nij-comp}
  Let $H$ and $H_1$ be nondegenerate Koszul-Vinberg structures on a left-symmetric algebroid $(A,\cdot_A,a_A)$. If $N=H^\sharp_1\circ (H^\sharp)^{-1}$ is a Nijenhuis operator on $A$, then $H$ and $H_1$ are compatible.
\end{thm}
\pf Since $N$ is a Nijenhuis operator on $A$, we have
$$NH^\sharp(\alpha)\cdot_A NH^\sharp(\beta )=N\big(NH^\sharp(\alpha)\cdot_A H^\sharp(\beta )+H^\sharp(\alpha)\cdot_A NH^\sharp(\beta )\big)-N^2(H^\sharp(\alpha)\cdot_A H^\sharp(\beta )).$$
By the assumption that $H$ and $H_1$ are Koszul-Vinberg structures,
\begin{eqnarray*}
H^\sharp(\alpha)\cdot_A H^\sharp(\beta)&=&H^\sharp(\huaL_{H^\sharp(\alpha)}\beta-R_{H^\sharp(\beta)}\alpha-\dM_A\langle H^\sharp(\alpha),\beta\rangle),\\
 NH^\sharp(\alpha)\cdot_A NH^\sharp(\beta) &=&NH^\sharp(\huaL_{NH^\sharp(\alpha)}\beta-R_{NH^\sharp(\beta)}\alpha-\dM_A\langle NH^\sharp(\alpha),\beta\rangle).
\end{eqnarray*}
Thus we obtain
\begin{eqnarray*}
  0&=&NH^\sharp(\alpha)\cdot_A NH^\sharp(\beta )-N\big(NH^\sharp(\alpha)\cdot_A H^\sharp(\beta )+H^\sharp(\alpha)\cdot_A NH^\sharp(\beta )\big)+N^2(H^\sharp(\alpha)\cdot_A H^\sharp(\beta ))\\
  &=&NH^\sharp(\huaL_{NH^\sharp(\alpha)}\beta-R_{NH^\sharp(\beta)}\alpha-\dM_A\langle NH^\sharp(\alpha),\beta\rangle)-N\big(NH^\sharp(\alpha)\cdot_A H^\sharp(\beta )\\
  &&+H^\sharp(\alpha)\cdot_A NH^\sharp(\beta )\big)+N^2H^\sharp(\huaL_{H^\sharp(\alpha)}\beta-R_{H^\sharp(\beta)}\alpha-\dM_A\langle H^\sharp(\alpha),\beta\rangle)\\
  &=&N\big(H^\sharp(\huaL_{NH^\sharp(\alpha)}\beta-R_{NH^\sharp(\beta)}\alpha-\dM_A\langle NH^\sharp(\alpha),\beta\rangle)+NH^\sharp(\huaL_{H^\sharp(\alpha)}\beta-R_{H^\sharp(\beta)}\alpha\\
  &&-\dM_A\langle H^\sharp(\alpha),\beta\rangle)-NH^\sharp(\alpha)\cdot_A H^\sharp(\beta )-H^\sharp(\alpha)\cdot_A NH^\sharp(\beta )\big).
\end{eqnarray*}
Because $N$ is nondegenerate,
\begin{eqnarray*}
&&NH^\sharp(\alpha)\cdot_A H^\sharp(\beta )+H^\sharp(\alpha)\cdot_A NH^\sharp(\beta )\\
 &=& H^\sharp(\huaL_{NH^\sharp(\alpha)}\beta-R_{NH^\sharp(\beta)}\alpha-\dM_A\langle NH^\sharp(\alpha),\beta\rangle)+NH^\sharp(\huaL_{H^\sharp(\alpha)}\beta-R_{H^\sharp(\beta)}\alpha-\dM_A\langle H^\sharp(\alpha),\beta\rangle),
\end{eqnarray*}
which implies that $H$ and $H_1$ are compatible.\qed\vspace{3mm}

The following proposition gives a local expression of Koszul-Vinberg structures, which will be used in the following examples.
\begin{pro}\label{pro:local coord}
Let $(M,\nabla)$ be a flat manifold and $\{x^1,\cdots,x^n\}$ be an affine coordinate
system with respect to $\nabla$. Then $H=\sum_{i,j=1}^nh_{i,j}\frac{\partial}{\partial x^i}\otimes \frac{\partial}{\partial x^j}\in\Sym^2(TM)$ is a Koszul-Vinberg structure on $M$ if and only if $$\sum_{l=1}^n(h_{jl}\frac{\partial h_{ik}}{\partial x^l}-h_{il}\frac{\partial h_{jk}}{\partial x^l})=0.$$
\end{pro}

\begin{ex}\label{ex:KVN-examp}{\rm
  For $\Real^2$ with the affine coordinate system $\{x,y\}$. Then, by a direct calculation in coordinates, $H$ and $H_1$ given by respectively
  $$H=\frac{\partial}{\partial x}\otimes \frac{\partial}{\partial x}+\frac{\partial}{\partial y}\otimes \frac{\partial}{\partial y}$$
  and
  $$H_1=\frac{x^2+y^2}{2}\frac{\partial}{\partial x}\otimes \frac{\partial}{\partial x}+xy\frac{\partial}{\partial x}\otimes \frac{\partial}{\partial y}++xy\frac{\partial}{\partial y}\otimes \frac{\partial}{\partial x}+\frac{x^2+y^2}{2}\frac{\partial}{\partial y}\otimes \frac{\partial}{\partial y}$$
  are Koszul-Vinberg structures and they are compatible. It is obvious that $H$ is nondegenerate. By Theorem {\rm\ref{thm:Comp-GHN}}, $N=H_1^\sharp\circ (H^\sharp)^{-1}$ given by
  $$
N=\begin{bmatrix}\frac{x^2+y^2}{2}& xy\\ xy&\frac{x^2+y^2}{2}\end{bmatrix}
$$
is a Nijenhuis operator on the left-symmetric algebroid $T_\nabla \Real^2 $. Furthermore, $(H,N)$ and $(H_1,N)$ are $\KVN$-structures on $\Real^2$.}
\end{ex}
\begin{ex}{\rm
  Consider the domain $\Omega=\{(x,y)\in\Real^2\mid x>0,y>0\}$ with the affine coordinate system $\{x,y\}$. Then the Koszul-Vinberg structures $H$ and $H_1$ given by respectively
  $$H=x\frac{\partial}{\partial x}\otimes \frac{\partial}{\partial x}+y\frac{\partial}{\partial y}\otimes \frac{\partial}{\partial y}$$
  and
  $$H_1=\frac{x^2+y^2}{2}\frac{\partial}{\partial x}\otimes \frac{\partial}{\partial x}+xy\frac{\partial}{\partial x}\otimes \frac{\partial}{\partial y}+xy\frac{\partial}{\partial y}\otimes \frac{\partial}{\partial x}+\frac{x^2+y^2}{2}\frac{\partial}{\partial y}\otimes \frac{\partial}{\partial y}$$
  are nondegenerate on $\Omega$. It is not hard to check that  $H$ and $H_1$ are not compatible. By Theorem \ref{thm:Nij-comp}, $N=H_1^\sharp\circ (H^\sharp)^{-1}$ given by
  $$
N=\begin{bmatrix}\frac{x^2+y^2}{2x}& x\\ y&\frac{x^2+y^2}{2y}\end{bmatrix}
$$
is not a Nijenhuis operator on the left-symmetric algebroid $T_\nabla \Omega $. }
\end{ex}

\begin{pro}\label{pro:GNN-nondegenerate}
 Let $(H,N)$ be a $\KVN$-structure on a left-symmetric algebroid $(A,\cdot_A,a_A)$. If $H$ is nondegenerate, then
 \begin{itemize}
\item[$\rm(i)$]$(H,N^k)$ is a $\KVN$-structure;
\item[$\rm(ii)$]$(H_{N^{k}},N^k)$ is a $\KVN$-structure, where $(H_{N^{k}})^\sharp=N^k\circ H^\sharp$.
\end{itemize}
\end{pro}
\pf Since $(H,N)$ is a $\KVN$-structure on the left-symmetric algebroid $(A,\cdot_A,a_A)$, by Theorem \ref{thm:hierarchy}, $H$ and $H_{N^{k}}$ are compatible $\KVN$-structures. Then by the condition that $H$ is nondegenerate and Theorem \ref{thm:Comp-GHN}, the conclusions follow immediately.\qed

\section{$\KVB$-structures and complementary symmetric $2$-tensors }\label{sec:MC}
In this section, we introduce the notion of a $\KVB$-structure on a left-symmetric algebroid, which consists of a Koszul-Vinberg structure and a symmetric $2$-cocycle such that some compatibility conditions hold. For any $\frkB\in\Sym^2(A^*)$, the bundle map  $B^\natural:A\longrightarrow A^*$  is defined by $\frkB^\natural(x)(y)=\frkB(x,y)$ for all $x,y\in\Gamma(A)$.
\begin{defi}
Let $H\in\Sym^2(A)$ be a Koszul-Vinberg structure and $\frkB\in\Sym^2(A^*)$ a   $2$-cocycle on a left-symmetric algebroid $(A,\cdot_A,a_A)$.
  Then $(H,\frkB)$ is called a {\bf $\KVB$-structure} if $\frkB_N $ is also a  $2$-cocycle, where $N=H^\sharp\circ \frkB^\natural$ and $\frkB_N\in\Sym^2(A^*)$ is defined by $\frkB_N(x,y)=\frkB( Nx,y)$ for all $x,y\in\Gamma(A)$.
\end{defi}
By the proof of Proposition 2.1 in \cite{Geometry of Hessian structures}, we have
\begin{pro}
  Let $(M,\nabla)$ be a flat manifold with $\{x^1,\cdots,x^n\}$ an affine coordinate
system with respect to $\nabla$ and $\frkB\in\Sym^2(T^*M)$. Then the following conditions are equivalent:
\begin{itemize}
\item[$\rm(i)$]$\frkB$ is a  $2$-cocycle on the left-symmetric algebroid $ T_\nabla M$;
\item[$\rm(ii)$]$\frac{\partial \frkB_{ij}}{\partial x^k}=\frac{\partial \frkB_{kj}}{\partial x^i}$;
\item[$\rm(iii)$]There exists a function $\varphi\in \CWM$ such that $\frkB$ can be locally expressed by
$\frkB_{ij}=\frac{\partial^2\varphi}{\partial x^i\partial x^j}$.
\end{itemize}
\end{pro}
\begin{ex}{\rm For $\Real^3$ with the affine coordinate system $\{x,y,z\}$. Define $H$ and $\frkB$ by
  $$H=\frac{\partial}{\partial x}\otimes \frac{\partial}{\partial x}+\frac{\partial}{\partial y}\otimes \frac{\partial}{\partial y},\quad \frkB=xdx\otimes dx+ydy\otimes dy+zdz\otimes dz.$$
It is easy to check that $H$ is a Koszul-Vinberg structure and $\frkB$ is a symmetric $2$-cocycle on  $T_\nabla\Real^3$ with the function $\varphi=\frac{1}{6}(x^3+y^3+z^3)$. Furthermore, $\frkB_N $ given by
$$\frkB_N=x^2dx\otimes dx+y^2dy\otimes dy$$
 is also a  $2$-cocycle  on $T_\nabla\Real^3$ with the function $\varphi=\frac{1}{12}(x^4+y^4)$. Thus $(H,\frkB)$ is a $\KVB$-structure.
}
\end{ex}
\begin{defi}
 Let $H\in\Sym^2(A)$ be a Koszul-Vinberg structure on a left-symmetric algebroid $(A,\cdot_A,a_A)$. A $2$-tensor $\frkB\in\Sym^2(A^*)$ is called a {\bf complementary symmetric $2$-tensor} for $H$ if
 $$
 \llbracket\frkB,\frkB\rrbracket_{A^*}=0,
 $$
 where $\llbracket-,-\rrbracket_{A^*}$ is given by \eqref{brac2} and the left-symmetric algebroid structure on $A^*$ induced by $H$ is given by Theorem \ref{thm:LSBi-H}.
\end{defi}

\emptycomment{
Similar to the definition of a strong Maurer-Cartan equation on a Lie bialgebroid introduced in \cite{lwx}, we give the following definition
\begin{defi}
Let $(A,A^*)$ be a left-symmetric bialgebroid and $\frkB\in\Sym^2(A^*)$.  The following equation
$$\delta_{A}\frkB=\llbracket \frkB,\frkB\rrbracket_{A^*}=0$$
 is called a {\bf strong Maurer-Cartan equation}.
\end{defi}
\begin{pro}\label{pro:comptoMC}
  Let $H\in\Sym^2(A)$ be a Koszul-Vinberg structure on a left-symmetric algebroid $(A,\cdot_A,a_A)$ and $\frkB\in\Sym^2(A^*)$. Then $\frkB$ is a $\delta_A$-closed complementary symmetric $2$-tensor for $H$ if and only if $\frkB$ is a solution of the strong Maurer-Cartan equation on the left-symmetric bialgebroid $(A,A^*)$ induced by $H$ given by Theorem \ref{thm:LSBi-H}.
\end{pro}
}

\begin{pro}
  Let $H\in\Sym^2(A)$ be a Koszul-Vinberg structure on a left-symmetric algebroid $(A,\cdot_A,a_A)$. Then $\frkB\in\Sym^2(A^*)$ is a complementary symmetric $2$-tensor for $H$ if and only if
  \begin{eqnarray}
 \nonumber &&a_A(H^\sharp \frkB^\natural(y))\frkB(x,z)-a_A( H^\sharp\frkB^\natural(x))\frkB(y,z)-a_A(x)\frkB(H^\sharp\frkB^\natural(y),z)\\
 \nonumber&&+a_A(y)\frkB(H^\sharp \frkB^\natural(x),z)-\frkB(H^\sharp \frkB^\natural(y)\cdot_Az,x) -\frkB(y\cdot_A H^\sharp\frkB^\natural(z),x)\\
 \nonumber&&+\frkB(H^\sharp \frkB^\natural(x)\cdot_Az,y)+\frkB(x\cdot_A H^\sharp \frkB^\natural(z),y)+\frkB([H^\sharp \frkB^\natural(x),y]_A,z)\\
&&-\frkB([H^\sharp \frkB^\natural(y),x]_A,z)=0\label{eq:com2-form}
\end{eqnarray}
 for all~$x,y,z\in \Gamma(A).$
\end{pro}
\pf It follows from
\begin{eqnarray*}
\llbracket\frkB,\frkB\rrbracket_{A^*}(x,y,z)&=&a_A(H^\sharp \frkB^\natural(y))\frkB(x,z)-a_A( H^\sharp\frkB^\natural(x))\frkB(y,z)-a_A(x)\frkB(H^\sharp\frkB^\natural(y),z)\nonumber\\
 &&+a_A(y)\frkB(H^\sharp \frkB^\natural(x),z)-\frkB(H^\sharp \frkB^\natural(y)\cdot_Az,x) -\frkB(y\cdot_A H^\sharp\frkB^\natural(z),x)\nonumber\\
 &&+\frkB(H^\sharp \frkB^\natural(x)\cdot_Az,y)+\frkB(x\cdot_A H^\sharp \frkB^\natural(z),y)+\frkB([H^\sharp \frkB^\natural(x),y]_A,z)\\&&-\frkB([H^\sharp \frkB^\natural(y),x]_A,z)
\end{eqnarray*}
for all~$x,y,z\in \Gamma(A).$ \qed
\begin{thm}\label{thm:comptoMC}
  Let $H\in\Sym^2(A)$ be a Koszul-Vinberg structure on a left-symmetric algebroid $(A,\cdot_A,a_A)$ and $\frkB\in\Sym^2(A^*)$. Then $(H,\frkB)$ is a $\KVB$-structure  if and only if $\frkB$ is a  $2$-cocycle that complementary to $H$.
\end{thm}
\pf Assume that $(H,\frkB)$ is a $\KVB$-structure. Let $N=H^\sharp\circ \frkB^\natural$. Since $\frkB$ is $\delta_A$-closed, we have
\begin{eqnarray*}
  a_A(N(x))\frkB(y,z)-a_A(y)\frkB(N(x),z)-\frkB(y,N(x)\cdot_A z)+\frkB(N(x),y\cdot_A z)-\frkB([N(x),y]_A,z)&=&0;\\
   a_A(x)\frkB(N(y),z)-a_A(N(y))\frkB(x,z)-\frkB(N(y),x\cdot_A z)+\frkB(x,N(y)\cdot_A z)-\frkB([x,N(y)]_A,z)&=&0;\\
   a_A(x)\frkB(y,N(z))-a_A(y)\frkB(x,N(z))-\frkB(y,x\cdot_A N(z))+\frkB(x,y\cdot_A N(z))-\frkB([x,y]_A,N(z))&=&0,
\end{eqnarray*}
which implies that
\begin{eqnarray*}
 &&a_A(N(x))\frkB(y,z)+ a_A(x)\frkB(N(y),z)-a_A(N(y))\frkB(x,z)-a_A(y)\frkB(N(x),z)\\
 &&+a_A(x)\frkB(y,N(z))-a_A(y)\frkB(x,N(z))-\frkB(y,N(x)\cdot_A z)-\frkB(N(y),x\cdot_A z)\\
 &&-\frkB(y,x\cdot_A N(z))+\frkB(N(x),y\cdot_A z)+\frkB(x,N(y)\cdot_A z)+\frkB(x,y\cdot_A N(z))\\
 &&-\frkB([N(x),y]_A,z)-\frkB([x,N(y)]_A,z)-\frkB([x,y]_A,N(z))=0.
\end{eqnarray*}
Since $\frkB_N$ is also  $\delta_A$-closed, we have
$$ a_A(x)\frkB(N(y),z)-a_A(y)\frkB(N(x),z)-\frkB(N(y),x\cdot_A z)+\frkB(N(x),y\cdot_A z)-\frkB([x,y]_A,N(z))=0.$$
Thus
\begin{eqnarray*}
 &&a_A( N(x))\frkB(y,z)-a_A( N(y))\frkB(x,z)+a_A(x)\frkB(N(y),z)\\
 &&-a_A(y)\frkB(N(x),z)-\frkB(y, x\cdot_A N(z))+\frkB(x ,y\cdot_A N(z))\\
 &&-\frkB(y ,N(x)\cdot_Az)+\frkB(x, N(y)\cdot_Az)-\frkB(z,[N(x),y]_A)-\frkB(z,[x ,N(y)]_A)=0,
\end{eqnarray*}
which is \eqref{eq:com2-form} with $N=H^\sharp\circ \frkB^\natural$ and thus $\frkB$ is a complementary symmetric $2$-tensor for $H$.

The converse part can be proved similarly. We finish the proof.\qed\vspace{3mm}

In the following, we will study the relation between $\KVN$-structures and $\KVB$-structures. Now we give a useful lemma.
\begin{lem}
 Let $H$ be a Koszul-Vinberg structure on a left-symmetric algebroid $(A,\cdot_A,a_A)$. Then $\frkB\in\Sym^2(A^*)$ is a symmetric $2$-cocycle that complementary to $H$ if and only if
 \begin{eqnarray}
   \label{eq:MC5}\frkB^\natural(x\cdot_A y)&=&\huaL_x\frkB^\natural(y)-R_y\frkB^\natural(x)-\dM_A\frkB(x,y),\\
     \label{eq:MC6}\frkB^\natural(x)\cdot^{H^\sharp} \frkB^\natural(y)&=&\frkB^\natural\big(H^\sharp\frkB^\natural(x)\cdot_A y+x\cdot_A H^\sharp\frkB^\natural(y)-H^\sharp\frkB^\natural(x\cdot_A y)\big).
 \end{eqnarray}
\end{lem}
\pf Since $\frkB$ is $\delta_A$-closed, we have $\delta_A\frkB(x,z,y)=0$, which is equivalent to
\begin{eqnarray*}
 \langle\huaL_x\frkB^\natural(y)-R_y\frkB^\natural(x)-\dM_A\frkB(x,y)-\frkB^\natural(x\cdot_A y),z\rangle=0.
\end{eqnarray*}
This is just \eqref{eq:MC5}.

Since $\llbracket \frkB,\frkB\rrbracket_{A^*}=0$, we have
\begin{eqnarray}\label{eq:MC3}
 \frkB^\natural(x)\cdot^{H^\sharp} \frkB^\natural(y)=\frkB^\natural(\huaL_{\frkB^\natural(x)}y-R_{\frkB^\natural(y)}x-\dM_{A^*}\frkB(x,y)),\quad\forall~x,y\in\Gamma(A).
\end{eqnarray}
Also, we have
\begin{eqnarray*}
&&\langle\huaL_{\frkB^\natural(x)}y-R_{\frkB^\natural(y)}x-\dM_{A^*}\frkB(x,y),\alpha\rangle\\
&=&a_{A^*}(\frkB^\natural(x))\langle y,\alpha\rangle-\langle\frkB^\natural(x)\cdot^{H^\sharp}\alpha,y \rangle+\langle x,\alpha\cdot^{H^\sharp}\frkB^\natural(x)\rangle-a_{A^*}(\alpha)\frkB(x,y)\\
&=&a_{A}(H^\sharp(\alpha))\frkB(x,y)-a_A(x)\frkB(H^\sharp(\alpha),y)-\frkB(x,H^\sharp(\alpha)\cdot_A y)-\frkB([H^\sharp(\alpha),x]_A,y)\\
&&+\langle H^\sharp\frkB^\natural(x)\cdot_A y,\alpha\rangle+\langle x\cdot_A H^\sharp\frkB^\natural(y),\alpha\rangle\\
&=&-\langle H^\sharp\frkB^\natural(x\cdot_A y),\alpha\rangle+\langle H^\sharp\frkB^\natural(x)\cdot_A y,\alpha\rangle+\langle x\cdot_A H^\sharp\frkB^\natural(y),\alpha\rangle\\
&=&\langle H^\sharp\frkB^\natural(x)\cdot_A y+x\cdot_A H^\sharp\frkB^\natural(y)-H^\sharp\frkB^\natural(x\cdot_A y),\alpha\rangle,
\end{eqnarray*}
where in the third equality we use the fact that $\delta_A\frkB(x,H^\sharp(\alpha),y)=0$. Thus, \eqref{eq:MC3} implies \eqref{eq:MC6}.

The converse part can be proved similarly. We finish the proof.\qed\vspace{3mm}
\begin{pro}\label{pro:ConNij-maxCon}
If $(H,\frkB)$ is a $\KVB$-structure on a left-symmetric algebroid $(A,\cdot_A,a_A)$, then $(H,N=H^\sharp\circ \frkB^\natural)$ is a $\KVN$-structure.

Conversely, if $(H,N)$ is a $\KVN$-structure and $H$ is nondegenerate, then $(H,\frkB)$ is $\KVB$-structure, where $\frkB$ is given by $\frkB^\natural=(H^\sharp)^{-1}\circ N$.
\end{pro}
\pf Now we assume that $(H,\frkB)$ is a $\KVB$-structure on a left-symmetric algebroid $(A,\cdot_A,a_A)$. First, applying $H^\sharp$ to both hands of \eqref{eq:MC6}, we obtain
\begin{eqnarray}\label{eq:MC7}
  H^\sharp\frkB^\natural(x)\cdot_A  H^\sharp\frkB^\natural(y)= H^\sharp\frkB^\natural(H^\sharp\frkB^\natural(x)\cdot_A y+x\cdot_A H^\sharp\frkB^\natural(y)-H^\sharp\frkB^\natural(x\cdot_A y)),
\end{eqnarray}
which implies that $N=H^\sharp\circ \frkB^\natural$ is a Nijenhuis operator on the left-symmetric algebroid $A$.

It is obvious that $N\circ H^\sharp=H^\sharp\circ \frkB^\natural\circ H^\sharp=H^\sharp\circ N^*.$
By a direct calculation, we have
\begin{eqnarray*}
&&\alpha\star^{H^\sharp}\beta-\alpha\cdot_{N^*}^{H^\sharp}\beta\\
&=&N^*(\alpha\cdot^{H^\sharp}\beta)-\huaL_{H^\sharp(\alpha)}N^*(\beta)+R_{H^\sharp(\beta)}N^*(\alpha)+\dM_AH(N^*(\alpha),\beta)\\
&=&\frkB^\natural \circ H^\sharp(\alpha\cdot^{H^\sharp}\beta)-\huaL_{H^\sharp(\alpha)}\frkB^\natural  H^\sharp(\beta)+R_{H^\sharp(\beta)}\frkB^\natural H^\sharp(\alpha)+\dM_AH(\frkB^\natural  H^\sharp(\alpha),\beta)\\
&=&\frkB^\natural(H^\sharp(\alpha)\cdot_A H^\sharp(\beta))-\huaL_{H^\sharp(\alpha)}\frkB^\natural  H^\sharp(\beta)+R_{H^\sharp(\beta)}\frkB^\natural H^\sharp(\alpha)+\dM_A\frkB( H^\sharp(\alpha),H^\sharp(\beta)).
\end{eqnarray*}
Therefore, \eqref{eq:MC5} implies that
$$\alpha\star^{H^\sharp}\beta=\alpha\cdot_{N^*}^{H^\sharp}\beta.$$ Thus $(H,N=H^\sharp\circ \frkB^\natural)$ is a $\KVN$-structure.

Conversely, since $N=H^\sharp\circ \frkB^\natural$ is a Nijenhuis operator on the left-symmetric algebroid $(A,\cdot_A,a_A)$, \eqref{eq:MC7} holds and
\begin{eqnarray*}
H^\sharp\big(\frkB^\natural(x)\cdot^{H^\sharp} \frkB^\natural(y)\big)&=&H^\sharp\frkB^\natural(x)\cdot_A  H^\sharp\frkB^\natural(y)\\
&=& H^\sharp\frkB^\natural(H^\sharp\frkB^\natural(x)\cdot_A y+x\cdot_A H^\sharp\frkB^\natural(y)-H^\sharp\frkB^\natural(x\cdot_A y)).
\end{eqnarray*}
Because $H^\sharp$ is nondegenerate, this implies \eqref{eq:MC6}.

Since $\alpha\star^{H^\sharp}\beta=\alpha\cdot_{N^*}^{H}\beta$, we have
$$\frkB^\natural(H^\sharp(\alpha)\cdot_A H^\sharp(\beta))-\huaL_{H^\sharp(\alpha)}\frkB^\natural \circ H^\sharp(\beta)+R_{H^\sharp(\beta)}\frkB^\natural \circ H^\sharp(\alpha)+\dM_A\frkB( H^\sharp(\alpha),H^\sharp(\beta))=0.$$
Also, $H^\sharp$ is nondegenerate, this implies \eqref{eq:MC5}.

By   Theorem \ref{thm:comptoMC}, $(H,\frkB)$ is a $\KVB$-structure. We finish the proof.\qed

\begin{ex}{\rm
   Consider the $\KVB$-structure $(H,\frkB)$ on $\Real^3$ given by the Example \ref{ex:KVN-examp}. Then by Theorem \ref{pro:ConNij-maxCon}, $(H,N)$ given by
  $$H=\frac{\partial}{\partial x}\otimes \frac{\partial}{\partial x}+\frac{\partial}{\partial y}\otimes \frac{\partial}{\partial y},\quad N=\begin{bmatrix}x& 0&0\\ 0&y&0\\0&0&0\end{bmatrix}$$
  is a $\KVN$-structure on $\Real^3$.
   }
\end{ex}
\section{Pseudo-Hessian-Nijenhuis structures}\label{sec:HN}

\begin{defi}
Let $\frkB\in\Sym^2(A^*)$ be a pseudo-Hessian structure and $N$ a Nijenhuis operator on a left-symmetric algebroid $({A},\cdot_{A},a_A)$. Then $(\frkB,N)$ is called a {\bf pseudo-Hessian-Nijenhuis structure}(or {\bf $\HN$-structure}) on the left-symmetric algebroid  $A$ if
\begin{equation}\label{eq:Hess1}
\frkB(N(x),y)=\frkB(x,N(y)),\quad \forall~x,y\in{\Gamma(A)}
\end{equation}
and $ \delta_A\frkB_N=0$, where $\frkB_N\in\Sym^2(A^*)$ is defined by $\frkB_N(x,y)=\frkB(N(x),y)$.
\end{defi}

\begin{thm}\label{thm:GHN-HN}
Let $(\frkB,N)$ be a $\HN$-structure on a left-symmetric algebroid  $({A},\cdot_{A},a_A)$. Then $(H,N)$, where $H^\sharp=(\frkB^\natural)^{-1}$, is a $\KVN$-structure on the left-symmetric algebroid $A$.

Conversely, if $(H,N)$ is a $\KVN$-structure and $H$ is nondegenerate, then $(\frkB,N)$, where $\frkB^\natural=(H^\sharp)^{-1}$, is a $\HN$-structure.
\end{thm}
\pf By Proposition \ref{pro:equivelent}, $H$ given by $(H^\sharp)^{-1}=\frkB^\natural$ is a Koszul-Vinberg structure on $A$.

By \eqref{eq:Hess1}, \eqref{eq:BNS1} follows immediately.
Now we prove \eqref{eq:BNS2}. First, by $\delta_A\frkB_N(x,z,y)=0$, we have
\begin{eqnarray}\label{eq:PHN-G1}
  a_A(x)\frkB(N(z),y)-B(N(z),x\cdot_A y)=a_A(z)\frkB(N(x),y)-\frkB(N(x),z\cdot_A y)+\frkB(N[x,z]_A,y).
\end{eqnarray}
Similarly, by $\delta_A\frkB(x,N(z),y)=0$, we have
\begin{eqnarray}\label{eq:PHN-G2}
  a_A(x)\frkB(N(z),y)-B(N(z),x\cdot_A y)=a_A(N(z))\frkB(x,y)-\frkB(x,N(z)\cdot_A y)+\frkB([x,N(z)]_A,y).
\end{eqnarray}
By \eqref{eq:PHN-G2}, \eqref{eq:PHN-G1} is equivalent to
\begin{eqnarray}
  &&a_A(z)\frkB(N(x),y)-\frkB(N(x),z\cdot_A y)+\frkB(N[x,z]_A,y)-a_A(N(z))\frkB(x,y)\nonumber\\
  &&+\frkB(x,N(z)\cdot_A y)-\frkB([x,N(z)]_A,y)=0.\label{eq:PHN-G3}
\end{eqnarray}
Let $\alpha=\frkB^\natural(x),\beta=\frkB^\natural(y)$, by \eqref{eq:PHN-G3},  we have
\begin{eqnarray*}
 &&\langle\alpha\star^{H^\sharp}\beta-\alpha\cdot_{N^*}^{H^\sharp}\beta,z\rangle\\
 &= &\langle N^*\huaL_{x}\frkB^\natural(y)-\huaL_{x}N^*\frkB^\natural(y)-N^*R_{y}\frkB^\natural(x)+R_{y}N^*\frkB^\natural(x)-N^*\dM_A\frkB^\natural(x,y)\\
 &&+\dM_A\frkB^\natural(N(x),y),z\rangle\\
 &=&a_A(z)\frkB(N(x),y)-\frkB(N(x),z\cdot_A y)+\frkB(N[x,z]_A,y)-a_A(N(z))\frkB(x,y)\\
 &&+\frkB(x,N(z)\cdot_A y)-\frkB([x,N(z)]_A)\\
 &=&0,
 \end{eqnarray*}
which implies that \eqref{eq:BNS2} holds for all $\alpha,\beta\in\Gamma(A^*)$. Thus, $(H,N)$ is a $\KVN$-structure.

The converse part can be proved similarly. We finish the proof.\qed

\begin{cor}
 Let $(\frkB,N)$ be a $\HN$-structure on a left-symmetric algebroid  $({A},\cdot_{A},a_A)$. Then $(\huaH,N)$, where $\huaH^\sharp=N\circ(\frkB^\natural)^{-1}$, is a $\KVN$-structure on the left-symmetric algebroid $A$.
\end{cor}
\pf Since $(\frkB,N)$ is a $\HN$-structure on the left-symmetric algebroid  $A$, then by Theorem \ref{thm:GHN-HN}, $(H,N)$,  where $H^\sharp=(\frkB^\natural)^{-1}$, is a $\KVN$-structure on the left-symmetric algebroid $A$. Furthermore, by Theorem \ref{thm:properties of GHN}, $H$ and $\huaH$  are compatible. Since $H$ is nondegenerate, by Theorem \ref{thm:Comp-GHN}, $(\huaH,N)$ is a $\KVN$-structure on the left-symmetric algebroid $A$. We finish the proof.\qed

\begin{ex}{\rm
  Consider the $\KVN$-structures $(H,N)$ and $(H_1,N)$ on $\Real^2$ given by the Example \ref{ex:KVN-examp}. Then by Theorem \ref{thm:GHN-HN},
  $(\frkB,N)$ given by
  $$\frkB=dx\otimes dx+dy\otimes dy,\quad N=\begin{bmatrix}\frac{x^2+y^2}{2}& xy\\ xy&\frac{x^2+y^2}{2}\end{bmatrix}$$
  is a $\HN$-structure on $\Real^2$. Since $H_1$ restricted to the domain $\Omega=\{(x,y)\in\Real^2\mid x^2-y^2>0\}$  is nondegenerate, by Theorem \ref{thm:GHN-HN},
  $(\frkB_1,N)$ given by
  $$\frkB_1=\frac{2x^2+2y^2}{(x^2-y^2)^2}dx\otimes dx-\frac{4xy}{(x^2-y^2)^2}dx\otimes dy-\frac{4xy}{(x^2-y^2)^2}dy\otimes dx+\frac{2x^2+2y^2}{(x^2-y^2)^2}dy\otimes dy$$
  and
  $$ N=\begin{bmatrix}\frac{x^2+y^2}{2}& xy\\ xy&\frac{x^2+y^2}{2}\end{bmatrix}$$
  is a $\HN$-structure on $\Omega$. Note that $\frkB_1$ can be locally expressed by the function $\varphi=\ln(x^2-y^2)$.
  }
\end{ex}

\begin{thm}\label{thm:HN-HB}
  If $(H,\frkB)$ is a $\KVB$-structure on a left-symmetric algebroid $(A,\cdot_A,a_A)$ and $\frkB$ is nondegenerate, then $(\frkB,N)$, where $N=H^\sharp\circ \frkB^\natural$, is a $\HN$-structure.

  Conversely, if $(\frkB,N)$ is a $\HN$-structure on a left-symmetric algebroid  $({A},\cdot_{A},a_A)$, then $(\huaH,\frkB)$, where $\huaH^\sharp=N\circ (\frkB^\natural)^{-1}$, is a $\KVB$-structure.
\end{thm}
\pf Since $(H,\frkB)$ is a $\KVB$-structure on the left-symmetric algebroid $(A,\cdot_A,a_A)$, $N=H^\sharp\circ \frkB^\natural$ is a Nijenhuis operator on $A$. By the definition of $\KVB$-structure, we have $\delta_A\frkB_N=0$. Thus $(\frkB,N)$ is a $\HN$-structure.

Conversely, since $(\frkB,N)$ is a $\HN$-structure on $A$, by Theorem \ref{thm:GHN-HN}, $(H,N)$, where $H^\sharp=(\frkB^\natural)^{-1}$, is a $\KVN$-structure. Thus $\huaH$ given by $\huaH^\sharp=N\circ (\frkB^\natural)^{-1}$ is a Koszul-Vinberg structure. It is obvious that $\delta_A\frkB_N=0$. Therefore, $(\huaH,\frkB)$ is a $\KVB$-structure. We finish the proof.\qed
\begin{pro}\label{pro:Herac-HN}
  Let $(\frkB,N)$ be a $\HN$-structure on a left-symmetric algebroid  $({A},\cdot_{A},a_A)$ . Then for all $k\in\Nat$, $\frkB_{N^k}$ defined by $(\frkB_{N^k})^\natural=\frkB^\natural\circ N^k$ are pseudo-Hessian structures and for any $k,l\in\Nat$, $\frkB_{N^k}$ and $\frkB_{N^l}$ are compatible in the sense that any linear combination of $\frkB_k$ and $\frkB_l$ are still  pseudo-Hessian structures.
\end{pro}
\pf By Theorem \ref{thm:GHN-HN}, $(H,N)$, where $H^\sharp=(\frkB^\natural)^{-1}$, is a $\KVN$-structure on the left-symmetric algebroid $A$. By the fact that $H$ is nondegenerate and Proposition \ref{pro:GNN-nondegenerate}, for any $k\in\Nat$, $(H,N^k)$ is a $\KVN$-structure. Also, by Theorem \ref{thm:GHN-HN}, $(\frkB,N^k)$ is a $\HN$-structure on $A$. Thus for any $k\in\Nat$, $\frkB_{N^k}$ is a pseudo-Hessian structure. The rest is direct. We finish the proof.\qed

\begin{cor}\label{cor:construction-HN}
 Let $N:A\rightarrow A$ be a bundle map and $\frkB$ a pseudo-Hessian structure on the left-symmetric algebroid $A$ such that \eqref{eq:Hess1} holds. Then $(\frkB,N)$ is a $\HN$-structure if and only if
 \begin{equation}
   \delta_A\frkB_N=\delta_A\frkB_{N^2}=0,
 \end{equation}
 where $\frkB_{N^i}(x,y)=\frkB(N^i(x),y)$ for $i=1,2$ and $x,y\in\Gamma(A)$.
\end{cor}
\pf If $(\frkB,N)$ is a $\HN$-structure, by Proposition \ref{pro:Herac-HN}, we have
\begin{equation*}
   \delta_A\frkB_N=\delta_A\frkB_{N^2}=0.
\end{equation*}
Conversely, let $H^\sharp=(\frkB^\natural)^{-1}$ and $\huaB=\frkB_N$ and note that $N=H^\sharp \circ \huaB^\natural$, then $\huaB_N=\frkB_{N^2}$. By hypothesis, $(H,\huaB)$ is a $\KVB$-structure and thus $N$ is a Nijenhuis operator on $A$. Therefore, $(\frkB,N)$ is a $\HN$-structure. We finish the proof.\qed\vspace{3mm}

We can relate the notion of a $\HN$-structure to that of a complementary symmetric $2$-tensor. The following proposition is a consequence of Theorem \ref{thm:comptoMC} and Theorem \ref{thm:HN-HB}.

\begin{pro}
  Let $({A},\cdot_{A},a_A)$ be a left-symmetric algebroid. If $\frkB\in\Sym^2(A^*)$ is a nondegenerate $2$-cocycle that complementary to a Koszul-Vinberg structure $H$, then $(\frkB,N)$, where $N=H^\sharp\circ \frkB^\natural$, is a $\HN$-structure.

  Conversely, if $(\frkB,N)$ is a $\HN$-structure, then $\huaH$ given by $\huaH^\sharp=N\circ (\frkB^\natural)^{-1}$ is a Koszul-Vinberg structure and $\frkB$ is a $2$-cocycle that complementary to $\huaH$.
\end{pro}
\begin{ex}{\rm
  Let $(\Real^n,\nabla)$ be the standard flat manifold and $\{x^1,\cdots,x^n\}$ be an affine coordinate
system with respect to $\nabla$. Let $\frkB\in\Sym^2(T^*\Real^n)$ be defined by
  $\frkB=\sum_{i=1}^n\dM x^i\otimes \dM x^i.$
  Then $\frkB$ is a pseudo-Hessian structure on $\Real^n$. Now consider the bundle map $N:TM\rightarrow TM$ by describing its conjugate operator $N^*:\Omega^1(\Real^n)\rightarrow \Omega^1(\Real^n)$ and put
  $$N^*\dM x^i=f_i(x^i)\dM x^i,\quad i=1,2,\cdots,n,$$
  where $f_i(x^i)$ is a smooth function on $\Real^n$ depending on the variable $x^i$.
  Note that $N$ is a Nijenhuis operator on $T_\nabla\Real^n$. In fact, it can be easily checked that
  $$\frkB(N(X),Y)=\frkB(X,N(Y)),\quad\forall~X,Y\in\frkX(\Real^n).$$
  For the formulas $\frkB_N(X,Y)=\frkB(N(X),Y)$ and $\frkB_{N^2}(X,Y)=\frkB(N^2(X),Y)$, we have
  \begin{eqnarray*}
    \frkB_N=\sum_{i=1}^nf_i(x^i)\dM x^i\otimes \dM x^i,\quad \frkB_{N^2}=\sum_{i=1}^n(f_i(x^i))^2\dM x^i\otimes \dM x^i,
  \end{eqnarray*}
  and therefore $\frkB_N$ and $\frkB_{N^2}$ are $2$-cocycles on $T_\nabla{\Real^n}$. Applying the result of Corollary \ref{cor:construction-HN}, we find that $(\frkB,N)$ is a $\HN$-structure.}
\end{ex}

{\footnotesize
\begin{thebibliography}{999}






\bibitem{Left-symmetric bialgebras}
C. Bai, Left-symmetric bialgebras and an analogue of the classical
Yang-Baxter equation. \emph{Commun. Contemp. Math.} 10 (2008),
221-260.

\bibitem{BBo}
S. Benayadi and M. Boucetta, On para-K${\rm\ddot{a}}$hler Lie algebroids and generalized pseudo-Hessian structures. \emph{Math. Nachr.}  292 (2019), 1418-1443.









\bibitem{Pre-lie algebra in geometry}
D. Burde, Left-symmetric algebras, or left-symmetric algebroids in geometry
and physics. \emph{Cent. Eur. J. Math.} 4 (2006), 323-357.













\bibitem{Dorfman1993}
I.  Dorfman, \emph{Dirac structures and integrability of nonlinear evolution equations,} Wiley, Chichester, 1993.

\bibitem{cohomology of pre-Lie}
A. Dzhumadil$'$daev, Cohomologies and deformations of right-symmetric
algebras. \emph{J. Math. Sci.} 93 (1999), 836-876.








\bibitem{G}
M. Gerstenhaber, The cohomology structure of an associative ring.
\emph{Ann. of Math.} 78 (1963), 267-288.


\bibitem{Kim}
H. Kim, Complete left-invariant affine structures on nilpotent Lie groups. \emph{J. Diff. Geom.}
24 (1986), 373-394.

\bibitem{Kosmann1}
Y. Kosmann-Schwarzbach and F. Magri, Poisson-Nijenhuis structures. \emph{Ann. Inst.
Henri Poincar$\acute{e}$ A} 53 (1990), 35-81.


\bibitem{Kosmann2}
Y. Kosmann-Schwarzbach and V. Rubtsov, Compatible structures on
Lie algebroids and Monge-Amp${\rm\grave{e}}$re operators. \emph{Acta. Appl. Math.} 109
(2010), no. 1, 101-135.

\bibitem{Koszul1}
 J. L. Koszul, Domaines born\'es homogenes et orbites de groupes de transformations affines. \emph{Bull. Soc. Math. France} 89 (1961), 515-533.







\bibitem{Lichnerowicz}
A. Lichnerowicz  and A. Medina, On Lie groups with left-invariant symplectic or K$\rm\ddot{a}$hlerian structures. \emph{Lett.
Math. Phys.} 16 (1988), 225-235.

\bibitem{LiuShengBaiChen}
J. Liu, Y. Sheng, C. Bai and Z. Chen, Left-symmetric algebroids.   \emph{Math. Nachr.}   289 (2016), No. 14-15, 1893-1908.

\bibitem{lsb}
J. Liu, Y. Sheng and C. Bai, Pre-symplectic algebroids and their applications.   \emph{Lett. Math. Phys.} 108 (2018), 779-804.

\bibitem{lsb2}
J. Liu, Y. Sheng and C. Bai, Left-symmetric bialgebroids and their corresponding Manin triples.  \emph{Diff. Geom. Appl.} 59 (2018), 91-111.

\bibitem{lwx}
Z. Liu, A. Weinstein and P. Xu. Manin triples for Lie bialgebroids.
 \emph{J. Diff. Geom.} 45 (1997), 547-574.

\bibitem{General theory of Lie groupoid and Lie algebroid}
K. C. H. Mackenzie, General theory of Lie groupoids and Lie
algebroids.  \emph{Lecture Note Series, $213$. London Mathematical
Society.}, Cambridge University Press, Cambridge, 2005.

\bibitem{Lie bialgebroid} K. C. H. Mackenzie and P. Xu, Lie bialgebroids and Poisson groupoids. \emph{Duke Math. J.} 73 (1994), 415-452.

\bibitem{MM}
F. Magri  and C. Morosi,  A Geometrical characterization of integrable Hamiltonian systems through
the theory of Poisson-Nijenhuis manifolds. \emph{Quaderno S19, University of Milan}, 1984.


\bibitem{MT}S. Majid and W. Tao, Noncommutative Differentials on Poisson-Lie groups and left-symmetric algebras.  \emph{Pacific J. Math.} 284 (2016), 213-256.


\bibitem{Med}
A. Medina, Flat left-invariant connections adapted to the automorphism structure
of a Lie group. \emph{J. Diff. Geom.} 16 (1981), 445-474.


\bibitem{Milnor}
J. Milnor, Curvatures of left invariant metrics on Lie groups. \emph{Adv. Math.} 21 (1976), 293-329.




\bibitem{Boyom1}
M. Nguiffo Boyom,  Cohomology of Koszul-Vinberg algebroids and Poisson manifolds. I. \emph{Banach Center Publ.} 54 (2001), 99-110.

\bibitem{Boyom2}
M. Nguiffo Boyom, KV-cohomology of Koszul-Vinberg algebroids and Poisson manifolds. \emph{Internat. J. Math.} 16 (2005), no. 9, 1033-1061.


\bibitem{NiBai}
X. Ni and C. Bai, Pseudo-Hessian Lie algebras and L-dendriform bialgebras. \emph{J. Algebra} 400 (2014), 273-289.





\bibitem{Shima}
H. Shima, Homogeneous Hessian manifolds. \emph{Ann. Inst. Fourier(Grenoble)} 30 (1980), 90-128.

\bibitem{Geometry of Hessian structures}
H. Shima, \emph{The geometry of Hessian structures.} World Scientific Publishing Co. Pte. Ltd., Hackensack, NJ, 2007. xiv+246 pp.






\bibitem{Vinb} E. B. Vinberg, Convex homogeneous cones. \emph{Transl. Moscou Math. Soc.} 12 (1963), 340-403.

\bibitem{Vaisman}
I. Vaisman, Complementary 2-forms of Poisson structures. \emph{Compos. Math.} 101
(1996), 55-75.

\bibitem{WBLS}
Q. Wang, Y. Sheng, C. Bai and J. Liu, Nijenhuis operators on pre-Lie algebras.  \emph{Commun. Contemp. Math.} doi: 10.1142/S0219199718500505.

\end{thebibliography}
}
\end{document}